\documentclass[12pt]{article}

\usepackage{amssymb, amsthm, amsgen, amscd, amssymb, amsmath, amsxtra}

\usepackage{graphics, epsfig, color}

\usepackage[matrix, arrow, curve]{xy}

\usepackage{amsmath,amscd}

\usepackage[pdftex, plainpages=false, pdfpagelabels]{hyperref}

\hypersetup{
    linktocpage=true,
    colorlinks=true,
    bookmarks=true,
    citecolor=blue,
    urlcolor=blue,
    linkcolor=blue,
    citebordercolor={1 0 0},
    urlbordercolor={1 0 0},
    linkbordercolor={.7 .8 .8},
    breaklinks=true,
    pdfpagelabels=true,
    }

\newtheorem{defn}{Definition}

\newtheorem{thm}{Theorem}
\theoremstyle{remark}

\def\bee{\begin{equation*}}

\def\eee{\end{equation*}}

\def\be{\begin{equation}}

\def\ee{\end{equation}}

\def\C{{\mathbb C}}

\def\R{{\mathbb R}}

\def\CP{{{\mathbb C}P}}

\definecolor{Maroon}{rgb}{0.8, 0.0, 0.0}

\title{\vskip-1in Triangles, Rotation, a Theorem and the Jackpot}
\author{Dave Auckly}

\date{}

\begin{document}
\bibliographystyle{plain}
\maketitle
\vskip-.2in
 \centerline{Mathematics Department, Kansas State University,
Manhattan, KS 66506 }
\centerline{{\tt{dav@math.ksu.edu} }}

\section*{Introduction}\label{intro}
The Atiyah-Singer index theorem is a wonderful theorem that has many varied applications.
The goal of this paper is to introduce this theorem and some of its applications to a broad audience in time to celebrate its golden anniversary. The theorem was announced in a paper published in 1963 \cite{AS1}. There is now a vast literature related to this theorem and its generalizations.
We are going describe a result about triangles, discuss linear differential equations, discuss the difference between circulation and rotation, state a theorem (guess which one), and hit a jackpot.

\section*{Triangles}
Consider the area of a sector of a sphere of radius $R$ having angle
$\theta$ as in the figure below. The area of such a sector is $S_\theta(R)=2R^2\theta$.
\begin{figure}[!ht]  
\hskip30bp
\includegraphics[width=45mm]{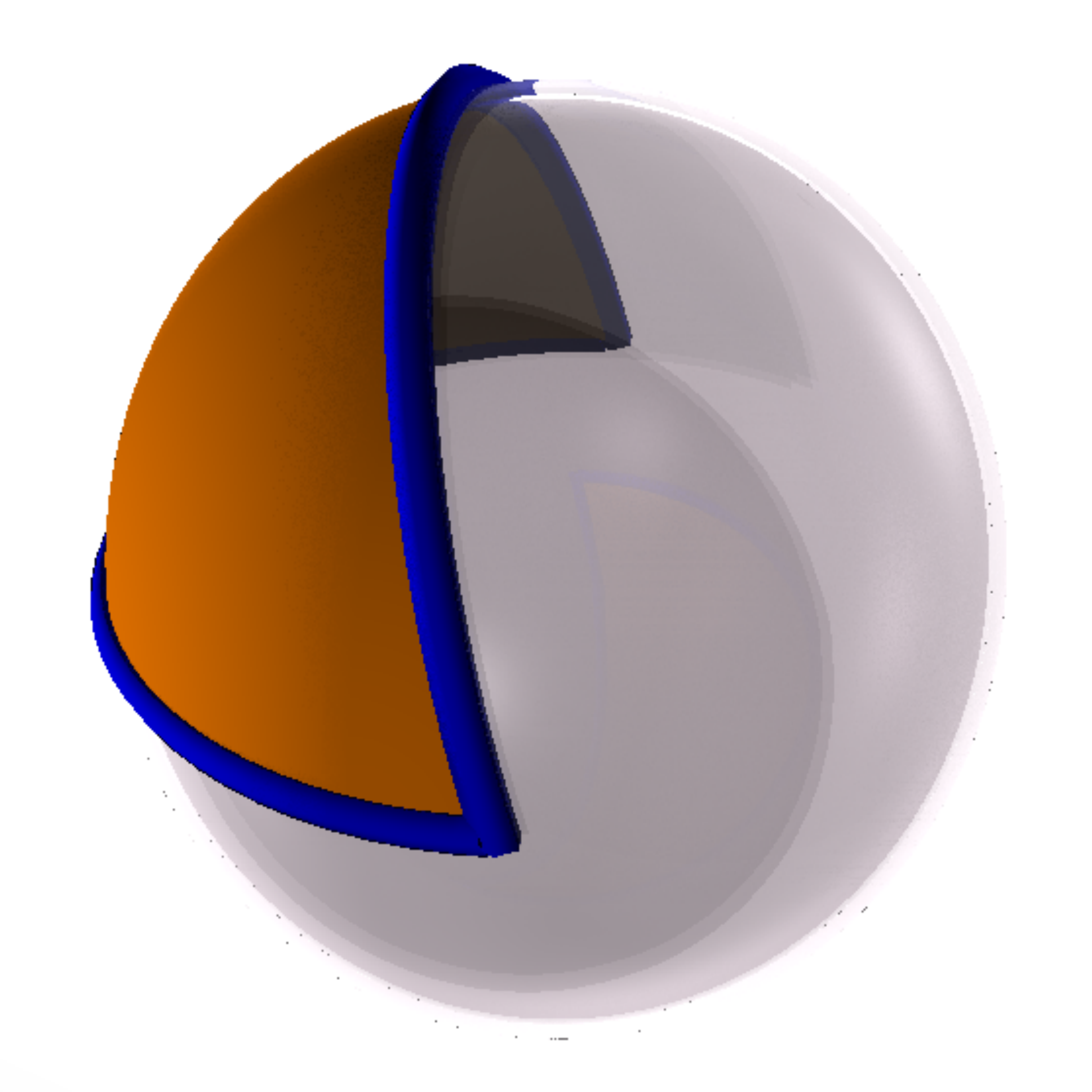}\hskip70bp\includegraphics[width=45mm]{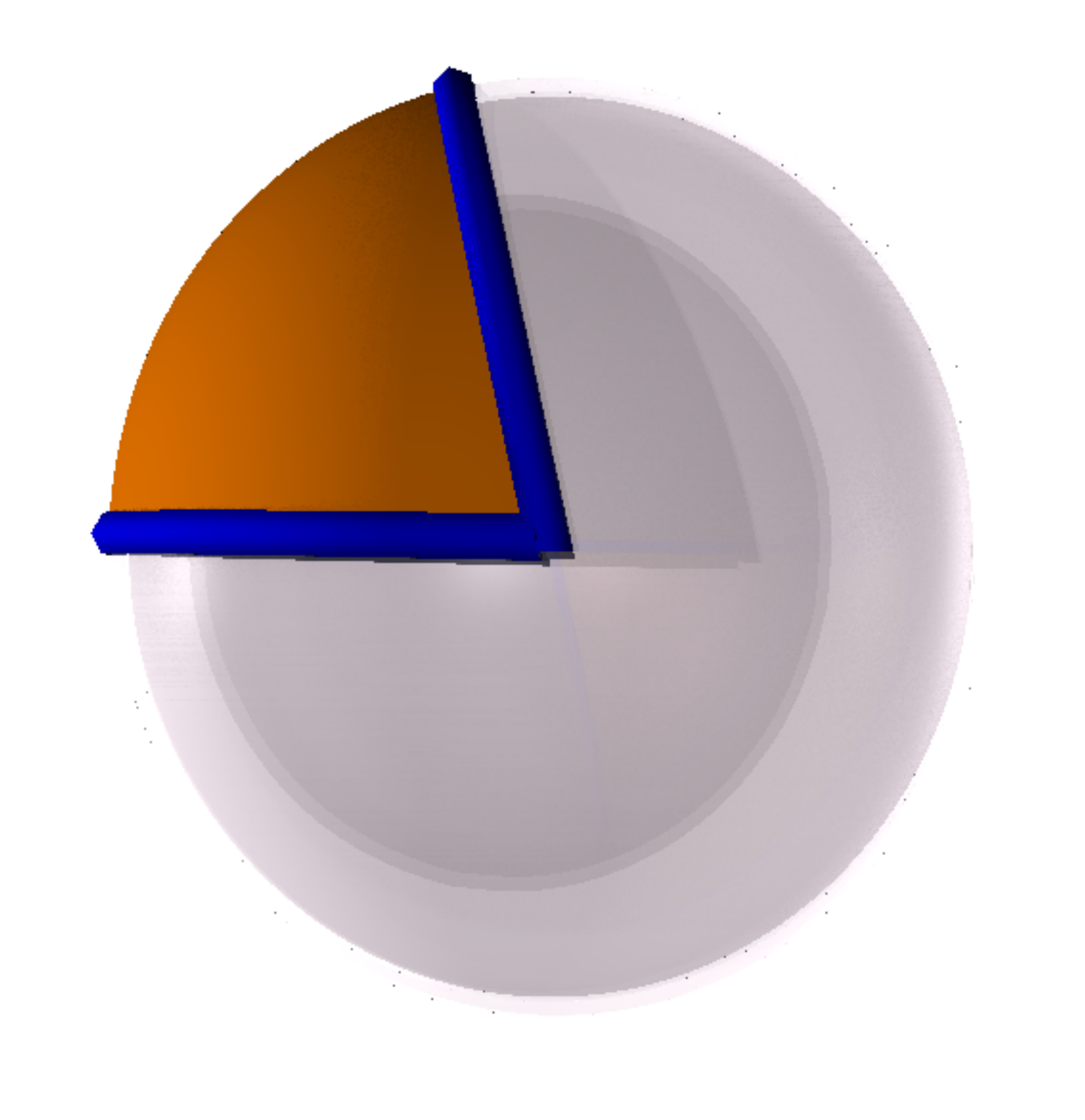}\caption{A sector}\label{sector}
\end{figure}

Now consider a spherical triangle on a sphere of radius $R$ as in the upper left of the following figure.
\begin{figure}[!ht]  
\[\begin{array}{cc}
\includegraphics[width=40mm]{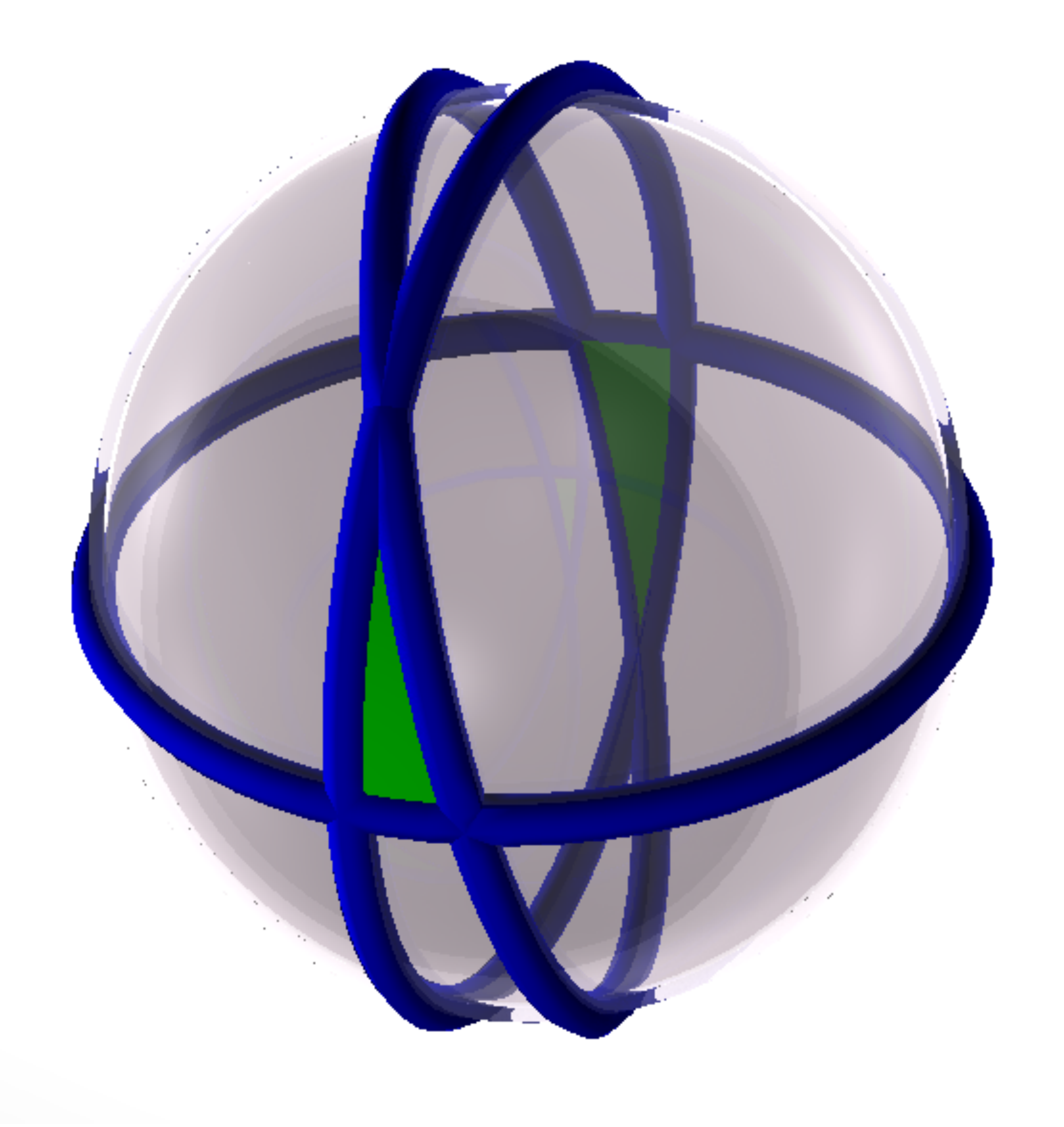} &\includegraphics[width=40mm]{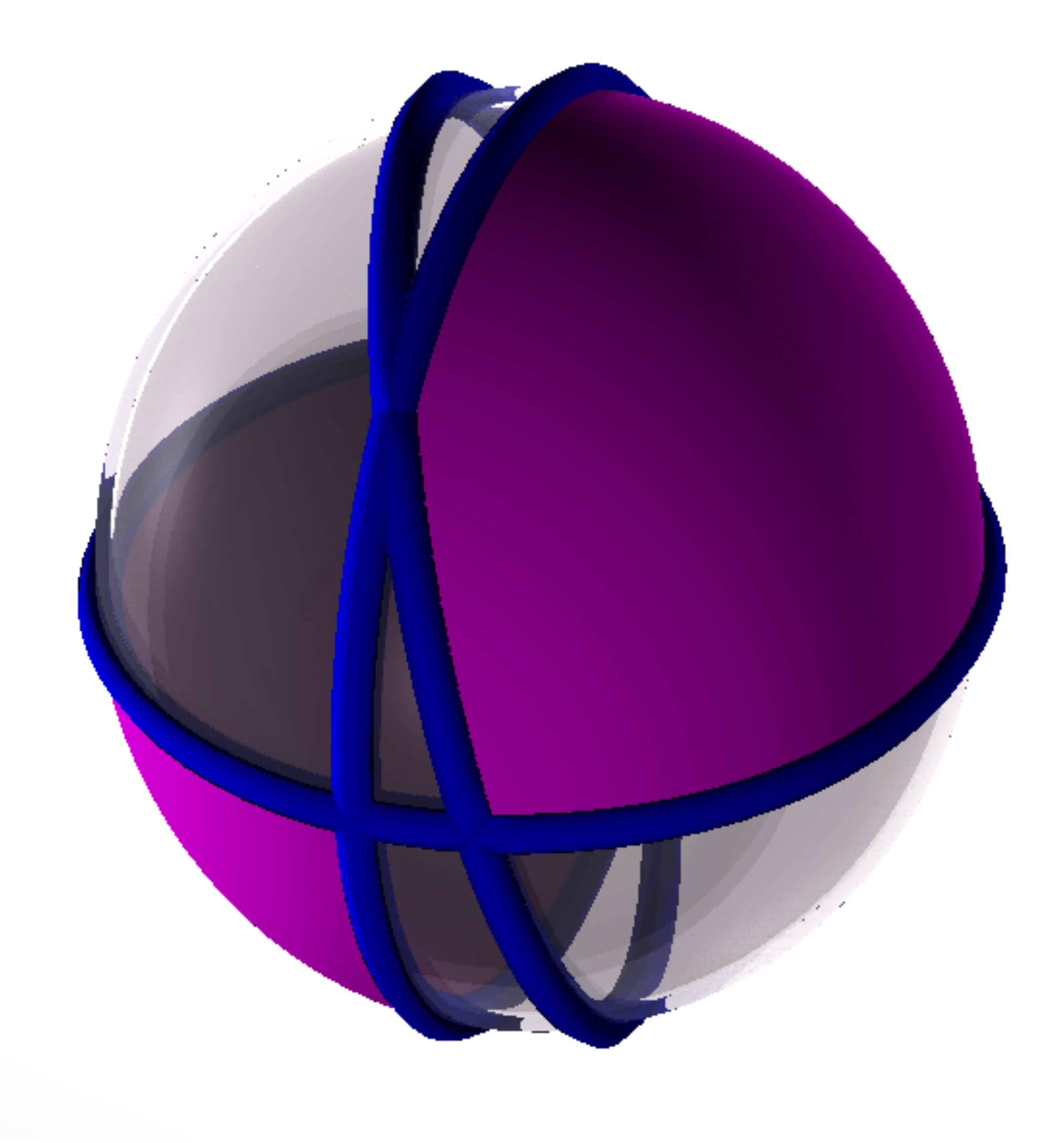}\\
\includegraphics[width=40mm]{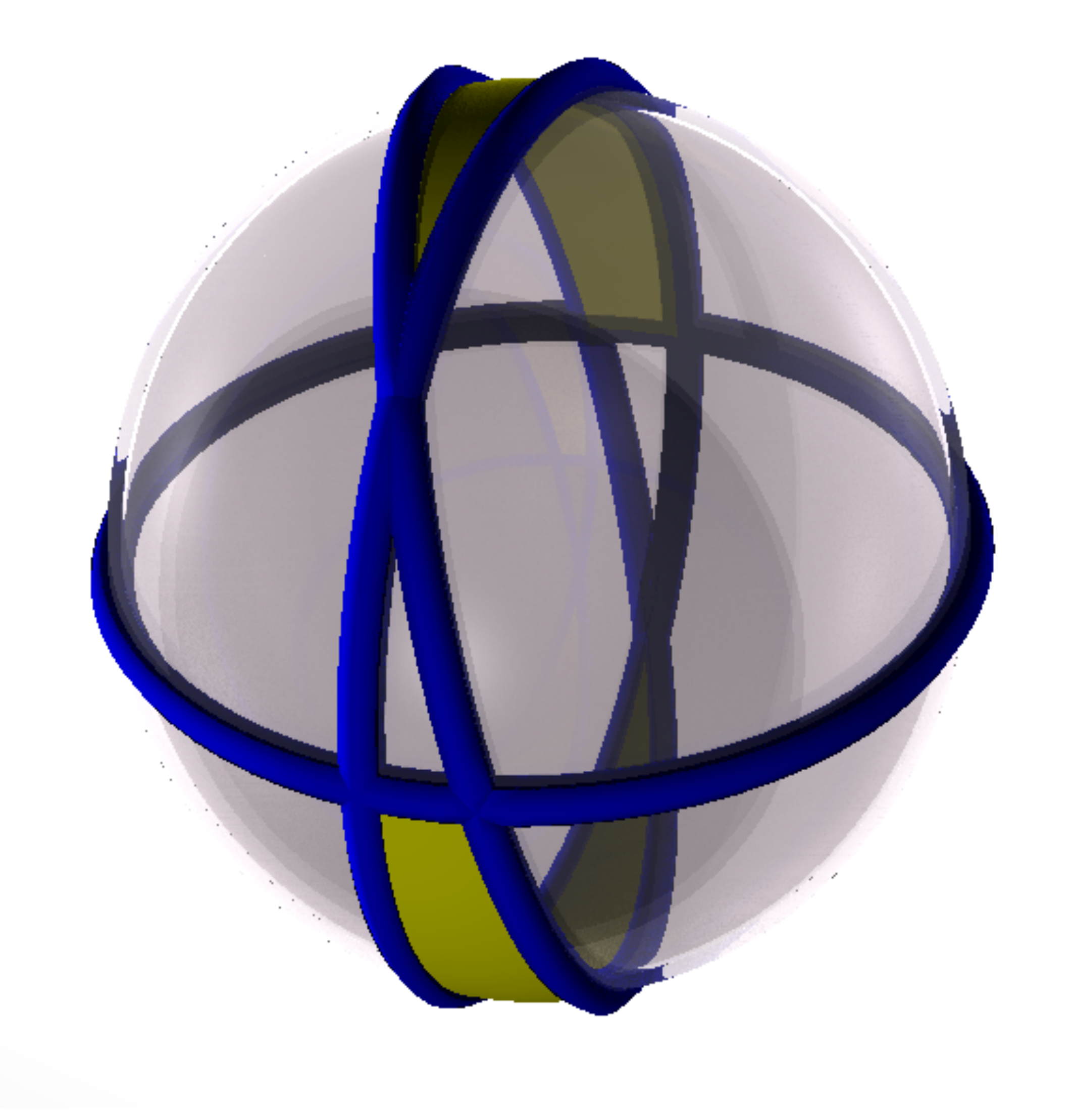} &\includegraphics[width=40mm]{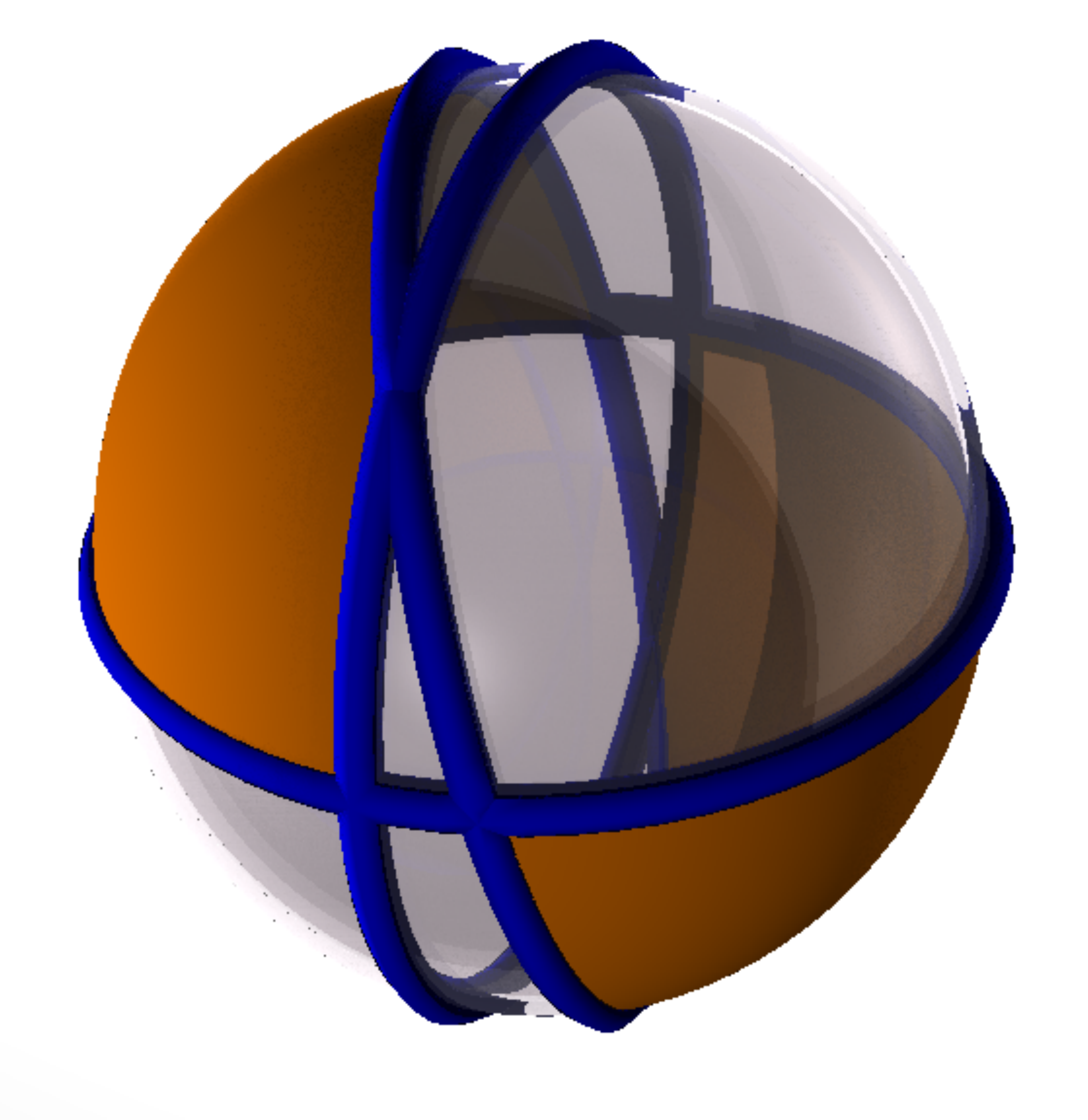}
\end{array} \begin{array}{c} \\ \includegraphics[width=50mm]{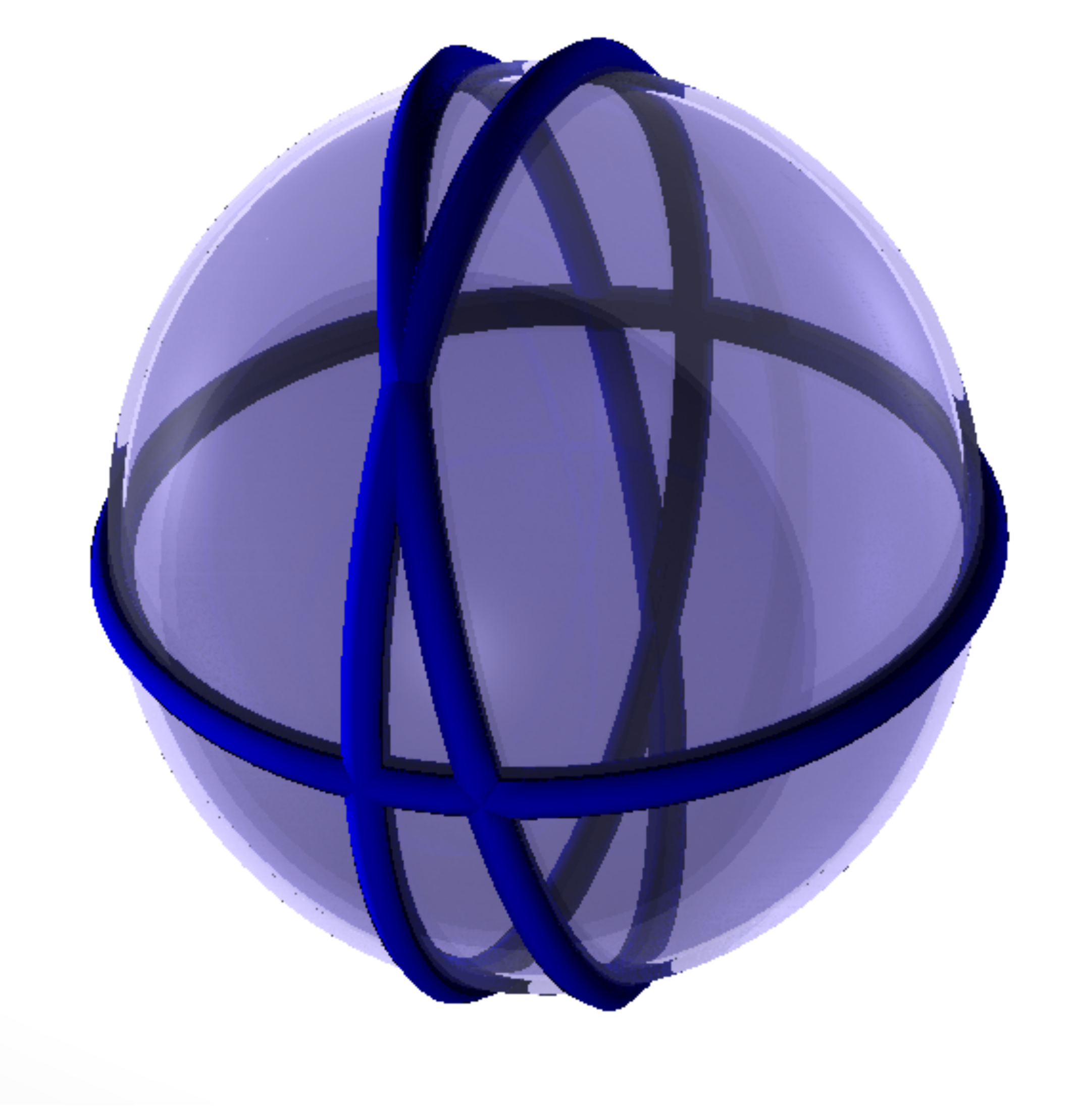} \\
\end{array}
\]
\caption{Spherical Triangle}\label{sphere}
\end{figure}
The great circles that make up the edges of the triangle divide the sphere into six truncated sectors and two copies of the triangle. Denote the angles of one of the two congruent spherical triangles by $\alpha$, $\beta$, and $\gamma$. Let $A$ denote the area of either triangle. The above figure proves:
\begin{equation*}
\begin{array}{lr}
2A&+2\left(S_\alpha(R)-A\right)\\
+2\left(S_\beta(R)-A\right)&+2\left(S_\gamma(R)-A\right)
\end{array}=4\pi R^2\,.
\end{equation*}
Substituting the formula for the area of a sector, combining like terms and dividing by $4\pi R^2$ gives:
\begin{equation*}
\alpha+\beta+\gamma-A\cdot\frac{1}{R^2}=\pi\,.
\end{equation*}
Notice that we can recover the formula for the sum of the angles in a triangle by constructing a suitable family of spherical triangles and taking the limit as $R\to 0$.

The fraction $\frac{1}{R^2}$ is known as the {\it sectional curvature}. It is denoted by $K$. With this notation the formula for the sum of the angles in a spherical triangle reads:
\begin{equation}\label{gb1}
\alpha+\beta+\gamma=\pi+A\cdot K\,.
\end{equation}

\begin{figure}[!ht]  
\hskip90bp
\includegraphics[width=75mm]{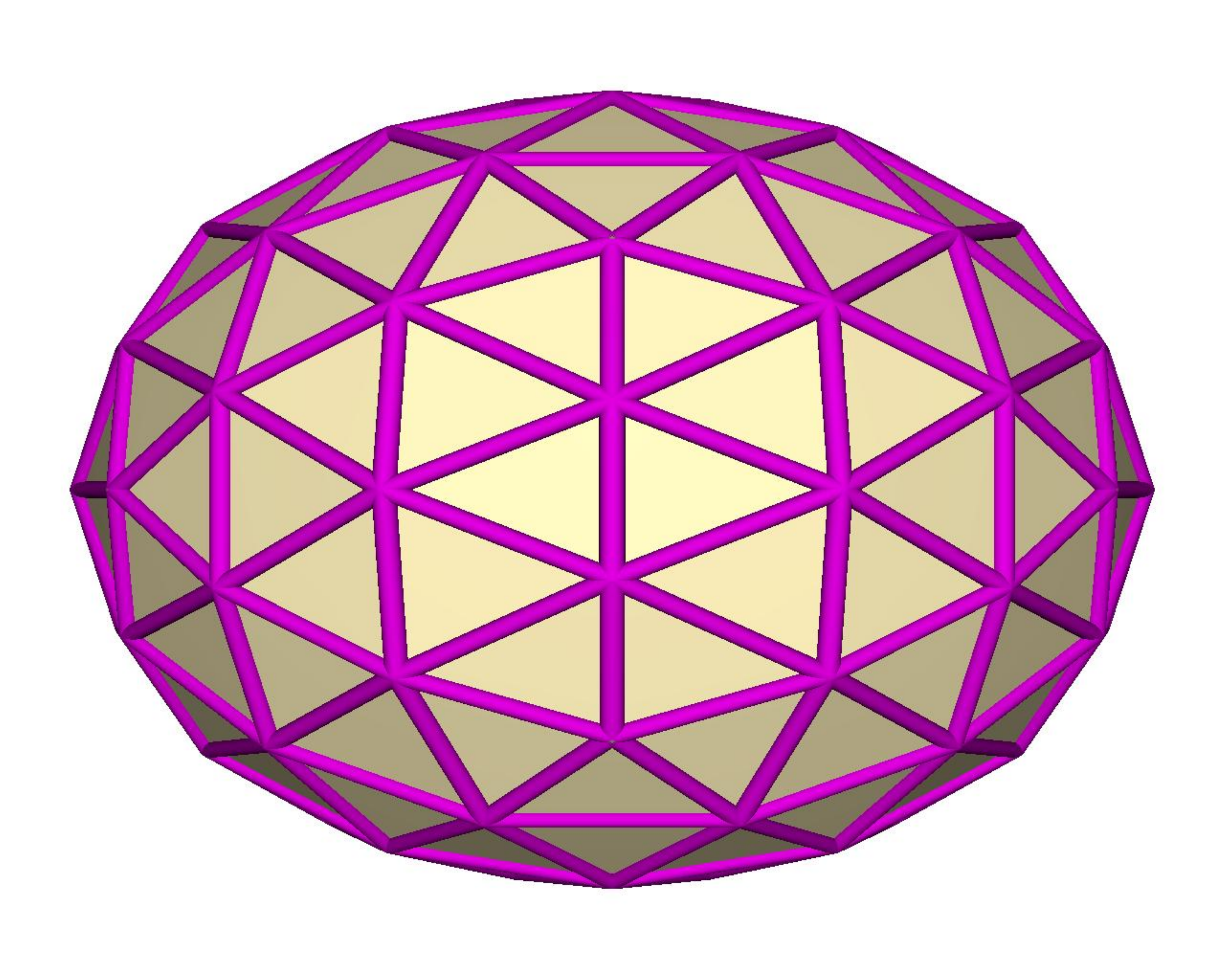}\caption{Triangulated surface}\label{triangulation}
\end{figure}

Now consider a triangulated surface as in figure \ref{triangulation}.
Summing the angle formula (\ref{gb1}) over all triangles leads to:
\begin{equation*}
\begin{array}{ccc}
\alpha_1 &+\beta_1 &+\gamma_1\\
+\alpha_2 &+\beta_2 &+\gamma_2\\
\cdot &\cdot &\cdot \\
+\alpha_n &+\beta_n &+\gamma_n\\
\end{array}
=
\begin{array}{c}
\pi+A_1K_1\\
+\pi+A_2K_2\\
\cdot\\
+\pi+A_nK_n\\
\end{array}
\end{equation*}
The sum of all of the angles that meet at a given vertex is $2\pi$. If we group all of the angles by vertex and let $V$ denote the number of vertices and $F$ the number of faces we obtain:
\begin{equation*}
\begin{aligned}
2\pi V&=\pi F+\frac{(A_1K_1+A_2K_2+\cdots+A_nK_3)}{(A_1+A_2+\cdots+A_n)}(A_1+A_2+\cdots+A_n)\\
&=\pi F+\overline{K}A=\pi F+\int_{\text{surface}}K\,d\,\text{Area}\,,
\end{aligned}
\end{equation*}
where
\[
\overline{K}=\frac{(A_1K_1+A_2K_2+\cdots+A_nK_3)}{(A_1+A_2+\cdots+A_n)}, \quad A=A_1+A_2+\cdots+A_n\,,
\]
are the average curvature and total area respectively.
This is a proof for constant positive sectional curvature and it certainly shows that the result is plausible in general.
Rearranging a bit gives:
\[
V-\frac12 F=\frac{1}{2\pi}\int_{\text{surface}}K\,d\,\text{Area}\,.
\]

There is something interesting about the previous formula. The left hand side does not depend on the curvature or geometry of the surface and the right hand side does not depend on the combinatorics of the triangulation. This means that the expression is independent of both. It is a topological invariant known as the Euler characteristic.

We can rewrite the left hand side a bit in order to make it manifestly invariant under
two basic moves that change the combinatorics. Imagine putting a red dot on each side of each edge. Thus there are $2E$ dots where $E$ denotes the number of edges. Since each face has three edges, there would be three dots in each face, so there are $3F$ dots and $3F=2E$. Adding and subtracting $E$ gives:
\[
V-\frac12 F=V-E+\frac32 F-\frac12 F=(V+F)-E\,.
\]
The quantity $\chi(\text{surface}):=(V+F)-E$ is known as the Euler characteristic.
Written in this way one can see that it will stay the same when an edge is subdivided into two by adding a new vertex. (The faces will no longer all be triangles.) It will also stay the same if a face is divided into two by adding a new edge connecting two vertices. One can use these moves to obtain an alternate proof that the Euler characteristic is an invariant.

To summarize we have given a plausibility argument for the Gauss-Bonett theorem:
\begin{equation}\label{gb2}
\chi(\text{surface})=\frac{1}{2\pi}\int_{\text{surface}}K\,d\,\text{Area}\,.
\end{equation}
The Gauss-Bonett theorem is good motivation for the Atiyah-Singer index theorem and follows from it as a corollary. We will outline a proof of it later.

\section*{Differential Equations}
The Gauss-Bonett theorem actually states something about differential equations, so we now turn to a discussion of linear PDEs. The existence and uniqueness of solutions of differential equations are two of the first issues addressed in the theory. Consider the equation $\text{grad}(u)={\bf V}$ as a linear partial differential equation for $u$. Since it is linear, solutions will be unique provided the only solution to $\text{grad}(u)=0$ is $u=0$. Of course any constant function has trivial gradient, so without imposing any boundary conditions solutions to our equation would not be unique. In fact there could even be more functions with trivial gradient depending on the topology of the domain. To understand this, let  $C^\infty_N(\Sigma)$ denote the space of smooth functions satisfying Neumann boundary conditions on a domain $\Sigma$, i.e. functions $w$ so that  $\hat n\cdot\text{grad}(w)=0$. Here $\hat n$ is the outer unit normal vector to the boundary. Set
\[
H^0(\Sigma):=\{w\in C^\infty_N(\Sigma) \, |\, \text{grad}(w)=0\}.
\]
If $\Sigma$ is the union of two disjoint intervals in $\R^1$, say $[0,1]$ and $[2,3]$, then $H^0(\Sigma)$ would be isomorphic to $\R^2$ with the isomorphism given by taking the ordered pair consisting the values of the function on each interval. In general $\text{dim}(H^0(\Sigma))$ is just the number of components of $\Sigma$.

Now turn to the existence of solutions of $\text{grad}(u)={\bf V}$. In this case one has to understand what are known as compatibility conditions. Recall that the rotation (or curl) of a vector field
\[
{\bf V}=f{\bf i}+g{\bf j}+h{\bf k}
\]
is given by
\[
\text{curl}({\bf V}):=\text{det}\left(
\begin{array}{ccc}
{\bf i}& {\bf j}& {\bf k}\\
\frac{\partial}{\partial x} &\frac{\partial}{\partial y} &\frac{\partial}{\partial z}\\
f& g& h
\end{array}\right)\,.
\]
Vector fields with vanishing curl are called {\it irrotational}. One computes that the curl of a gradient is trivial, thus one can only solve $\text{grad}(u)={\bf V}$ provided that the compatibility condition $\text{curl}({\bf V})=0$ is satisfied. This leaves us to solve the compatibility equation or more generally $\text{curl}({\bf V})=W$, as well as to understand how the original problem may be solved when the compatibility condition is satisfied.

Solving $\text{curl}({\bf V})={\bf W}$ appears to be analogous to our original problem solving
$\text{grad}(u)={\bf V}$, and indeed it is. One can see that the corresponding compatibility conditions in three dimensions are given by $\text{div}({\bf W})=0$. We do not, however, need to follow this thought any further.

The second question of solving the original gradient problem is more interesting to us. Since the gradient is essentially just a derivative, it is not surprising that a function with given gradient may be found by integration. The formula $u(\delta(1))-u(\delta(0))=\int_\delta \text{grad}(u)\cdot d\hat s$, where $d\hat s:={\bf i}\,dx+{\bf j}\,dy$ is the usual line element, and $\delta:[0,1]\to \Sigma$ will determine the value of the function $u$ in terms the value at any point the domain (assuming the region is path connected).  This meshes well with our analysis of the zeros of the gradient. To summarize the solution to $\text{grad}(u)={\bf V}$ is $u(x)=u(x_*)+\int_\delta V\cdot d\hat s$.  However, there is a problem -- different paths may lead to different values of the function $u$. This may be measured by the value of thee integral around a loop, consisting of taking one path to the point and a second path back.

The {\it circulation} around $\gamma$ is given by  $\int_\gamma {\bf V}\cdot d\hat s$.
The circulation of the gradient of any function about any closed curve is zero. In fact, any vector field with vanishing circulation is the gradient of a vector field as the formula in the previous paragraph shows. A non-constant periodic integral curve of a vector field is a good example of a curve with non-vanishing circulation. We see that a gradient field can have no such curves. This makes sense since it is impossible for one to come back to the starting point of a walk by always walking up.
\begin{figure}[!ht]  
\hskip110bp
\includegraphics[width=65mm]{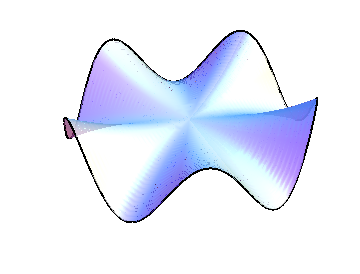}\caption{The gradient of a function has no circulation.}\label{gradient}
\end{figure}

\begin{figure}[!ht]  
\hskip110bp
\includegraphics[width=65mm]{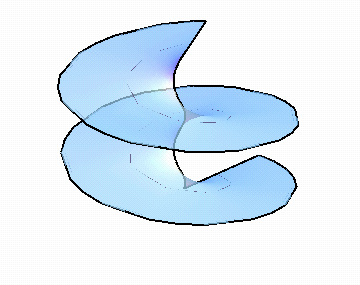}\caption{The ``gradient" of a helicoid  has non-zero circulation.}\label{helicoid}
\end{figure}

The concepts of rotation and circulation are related, but they are not the same. A vector field may be irrotational and still have non-trivial circulation. Even though such vector fields are not the gradient of any function, they may be viewed as the ``gradient" of a multi-valued function, such as  the helicoid in figure \ref{helicoid}. Here the walk up the ``spiral staircase edge" projects to a circle in the plane that is an integral curve of the ``gradient" of the helicoid. The circulation around this curve is non-trivial even though the vector field is irrotational.

Our discussion is about to take flight. It is reasonable to assume that airflow below one third of the speed of sound is irrotational. Thus for example, one curve may circulate around the wing of a plane and the other may wrap around the trailing vortex. As these two curves bound a common (bent) cylinder the corresponding circulations must be the same. The lift acting on the airplane is proportional to this circulation \cite{anderson}. Thus, one could estimate the lift of a plane based on a pair of timed photos taken from behind the plane. 
\begin{figure}[!ht]  
\hskip70bp
\includegraphics[width=95mm]{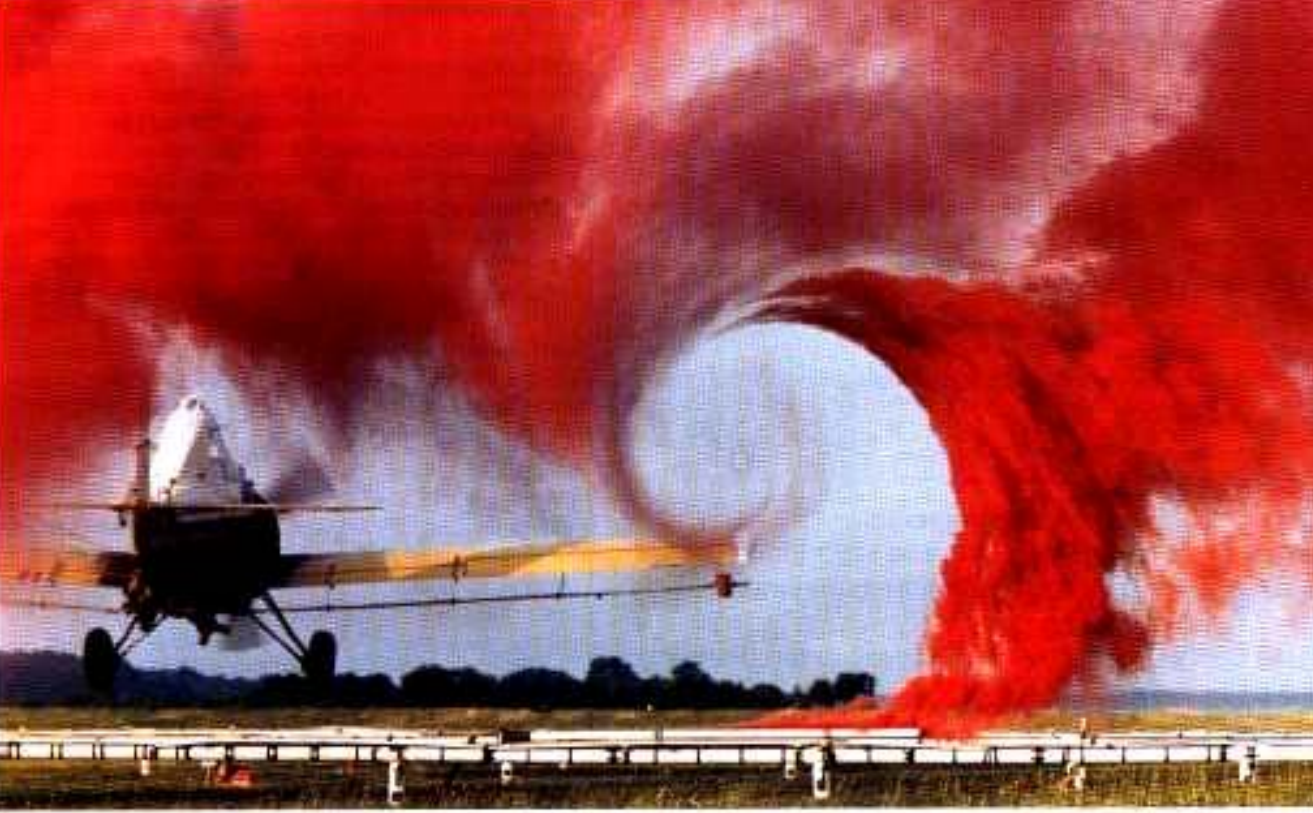}\caption{The circulation of the trailing vortex is proportional to the lift.}\label{plane}
\end{figure}

To see that the circulation of an irrotational vector field about two curves that cobound a surface must be the same, one applies Stokes' theorem:
\[
\int_{\gamma_1}{\bf V}\cdot d\hat s\,-\,\int_{\gamma_2}{\bf V}\cdot d\hat s=\int_{\partial\Sigma^\prime}{\bf V}\cdot d\hat s
=\int_{\Sigma^\prime}\text{curl}({\bf V})\cdot {\bf k}\,d\,\text{Area}=0\,.
\]

Let $\Sigma$ be a bounded $2$-dimensional domain with smooth boundary. Forcing the velocity field to be tangent to the boundary is natural in the context of fluid flows. Let $\mathcal{X}_{T}(\Sigma)$ denote the space of smooth vector fields on $\Sigma$ that satisfy ${\bf V}\cdot\hat n=0$. Since we are restricting to a $2$-dimensional domain, the curl is encoded in the one component ${\bf k}\cdot\text{curl}$. We can now analyze the equation $\text{curl}({\bf V})=0$ in much the same way we analyzed the equation $\text{grad}(u)=0$. Once again, the topology of the domain will determine the answer. Since the curl of the gradient is trivial, we are really interested in irrotational vector fields that are not gradient fields. The key is to construct the following quotient space that decares that the gradient of any function is trivial:
\[
H^1(\Sigma):=\{{\bf V} \in \mathcal{X}_{T}(\Sigma) \, | \, {\bf k}\cdot\text{curl}({\bf V}) =0\} /\text{Image}(\text{grad}:C^\infty_N(\Sigma)\to
\mathcal{X}_{T}(\Sigma))\,.
\]
This now shows why we chose Neumann boundary conditions for our original gradient equation. The gradient of a function in $C^\infty_N(\Sigma)$ is tangent to the boundary and thus satisfies the boundary conditions that are natural for fluid flow.

\begin{figure}[!ht]  
\hskip70bp
\includegraphics[width=85mm]{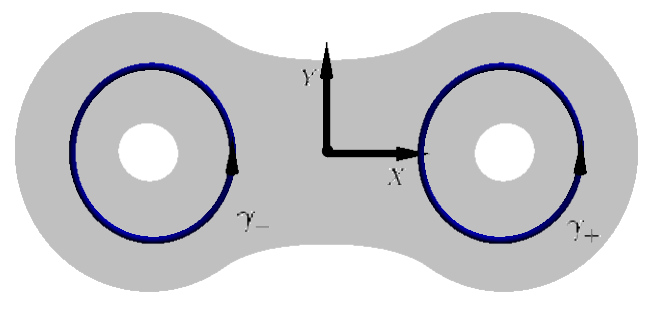}\caption{The domain $\Sigma$}\label{domain}
\end{figure}

It turns out that $\text{dim}(H^1(\Sigma))$ is the number of ``holes" in the domain $\Sigma$. As an example consider the domain $\Sigma$ displayed in figure \ref{domain}. The figure indicates a pair of simple closed curves in the domain. We will denote these curves by $\gamma_+$ and $\gamma_-$.
We will see that $H^1(\Sigma)$ is two dimensional in this example because $\Sigma$ has two ``holes," one around $(1,0)$ and one $(-1,0)$.  In fact the following function is an isomorphism.
\[
\Psi:H^1(\Sigma)\to \R^2; \quad \Psi({\bf V}):=(\int_{\gamma_-}{\bf V}\cdot d\hat{s},\int_{\gamma_+}{\bf V}\cdot d\hat{s})\,.
\]

To see that $\Psi$ is injective, assume that $\Psi({\bf V})=0$. It is  plausible (and in fact true) that any closed curve in $\Sigma$ together with some number of copies of $\gamma_-$ and $\gamma_+$ bounds a parameterized subsurface of $\Sigma$, say $\Sigma^\prime$. Stokes' theorem then implies that the circulation of any irrotational vector field about any curve in $\Sigma$ is zero. Thus the vector field ${\bf V}$ is the gradient of some function, and thus represents zero in $H^1(\Sigma)$.

The vector fields ${\bf V}_\pm:=\frac{1}{2\pi}\frac{(x\mp 1){\bf j}-y{\bf i}}{(x\mp 1)^2+y^2}$ may be used to show that $\Psi$ is surjective. They do not satisfy the proper boundary conditions, but this can be fixed by adding appropriate gradients.

The curl of any vector field $w={\bf k}\cdot\text{curl}({\bf V})$ satisfies Dirichlet boundary conditions on $\Sigma$, i.e. restrict to zero on the boundary for any ${\bf V}\in \mathcal{X}_{T}(\Sigma)$. The space of function satisfying Dirichlet boundary conditions on $\Sigma$ is denoted $C^\infty_D(\Sigma)$. To finish our analysis, we should address the ``two dimensional holes:"
\[
H^2(\Sigma):=C^\infty_D(\Sigma)/\text{Image}({\bf k}\cdot\text{curl}:\mathcal{X}_{T}(\Sigma)\to C^\infty_D(\Sigma))\,.
\]

To analyze this we will use the notion of the adjoint of a linear operator. The {\it adjoint} of a linear operator $D:E\to F$ between inner product spaces is the operator $D^*:F\to E$ such that $\langle Dx,y\rangle=\langle x,D^*y\rangle$. Clearly the adjoint of $D^*$ is $D$. The space $F/\text{Image}(D)$ is known as the cokernel of $D$. It is isomorphic to the kernel of $D^*$. This is especially easy to see in the case when the cokernel is trivial. A vector $y\in F$ is perpendicular to the image of $D$ if and only if $0=\langle Dx,y\rangle=\langle x,D^*y\rangle$ for every $x\in E$, and this is only true when $y$ is in the kernel of $D^*$. 

It is clear that $H^2(\Sigma)$ is the cokernel of ${\bf k}\cdot\text{curl}$. The pairing on $\mathcal{X}_{T}(\Sigma)$ is given by $\langle {\bf V},{\bf W} \rangle := \int_\Sigma {\bf V}\cdot {\bf W} \,d\,\text{Area}$. The pairing on $C^\infty_D(\Sigma)$ is given by $\langle u,w \rangle := \int_\Sigma uw \,d\,\text{Area}$. The following computation shows that the adjoint of ${\bf k}\cdot\text{curl}$ is $\text{grad}(\cdot)\times {\bf k}$. Since the only function with trivial gradient satisfying Dirichlet boundary conditions is zero, the cokernel of ${\bf k}\cdot\text{curl}$ is trivial.

The curl satisfies the following product rule:
\[
\text{curl}(w{\bf V})=w\,\text{curl}({\bf V})+\text{grad}(w)\times {\bf V}\,.
\]
For $w\in C^\infty_D(\Sigma)$ we have
\[
\begin{aligned}
0=\int_{\partial\Sigma}w\,{\bf V}\cdot d\hat s &= \int_\Sigma \text{curl}(w{\bf V})\cdot {\bf k}\,d\,\text{Area}\\
&=\int_\Sigma w\,\text{curl}({\bf V})\cdot {\bf k}\,d\,\text{Area}+
\int_\Sigma (\text{grad}(w)\times {\bf V})\cdot {\bf k}\,d\,\text{Area}\\
&=\int_\Sigma w\,\text{curl}({\bf V})\cdot {\bf k}\,d\,\text{Area}-
\int_\Sigma (\text{grad}(w)\times {\bf k})\cdot {\bf V}\,d\,\text{Area}\,.
\end{aligned}
\]

A complex of operators 
\[
\begin{CD}
E_0 @>D_0>> E_1 @>D_1>> E_2\,,
\end{CD}
\]
i.e. $D_1D_0=0$ may be collapsed to a single operator $D:E_0\oplus E_2 \to E_1$ given by $D=D_0+D_1^*$. One may see that $\text{ker}(D)\cong H^0\oplus H^2$ where $H^0=\text{ker}(D_0)$ and $H^2=\text{coker}(D_1)$. Indeed, $D_0(u)+D_1^*(w)=0$ implies $D_1D_1^*w=0$, so $\langle D_1^*w,D_1^*w\rangle=\langle D_1D_1^*w,w\rangle = 0$ giving $D_1^*w=0$ and therefore $D_0u=0$. Similarly, one may see that $\text{coker}(D)\cong H^1$, where $H^1=\text{ker}(D_1)/\text{image}(D_0)$.

Apply all of our discussion so far to the operator obtained by collapsing the $\text{grad}$-$\text{curl}$ complex,
\[
D:C^\infty_N(\Sigma)\oplus C^\infty_D(\Sigma)\to \mathcal{X}_{T}(\Sigma)\,,
\]
given by $D(u,w):=\text{grad}(u)+\text{grad}(w)\times {\bf k}$.
For a domain with $g$ holes we have
\begin{equation}\label{indd}
\begin{array}{l}\text{dim}(\text{ker}(D))=1\\
\text{dim}(\text{coker}(D))=g\,.
\end{array}
\end{equation}

\begin{figure}[!ht]  
\hskip70bp
\includegraphics[width=95mm]{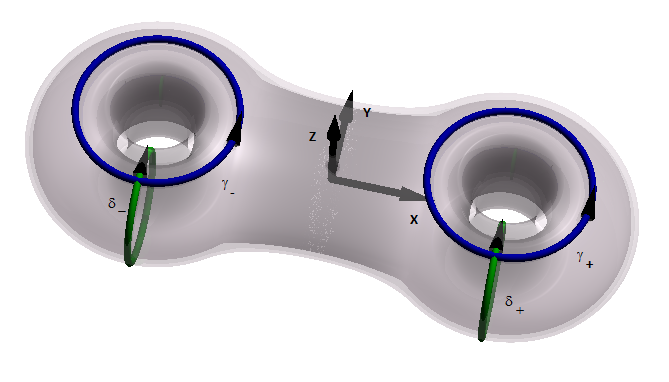}\caption{The surface $\widehat\Sigma$}\label{sighat}
\end{figure}

We wish to combine this with the Gauss-Bonnet theorem. Since the Gauss-Bonnet theorem applies to {\it closed} surfaces we need to construct a closed surface. We do so by doubling the surface $\Sigma$ into a surface $\widehat\Sigma$ as in figure \ref{sighat}.
We recover $\Sigma$ by projecting $\widehat\Sigma$ to the $xy$-plane. This projection identifies a point $(x,y,z)\in\widehat\Sigma$ with the point $\tau(x,y,z):=(x,y,-z)$. Smooth functions $u$ on $\widehat\Sigma$ such that $\tau^*u:=u\circ \tau = u$ correspond exactly to functions on $\Sigma$ satisfying Neumann boundary conditions i.e. $C^\infty_N(\Sigma)$. Similarly functions $w$ on $\widehat\Sigma$ such that $w\circ \tau = -w$ correspond exactly to functions on $\Sigma$ satisfying Dirichlet boundary conditions i.e. $C^\infty_D(\Sigma)$. It is also possible to identify the set of vector
fields on $\Sigma$ satisfying the given boundary conditions with $\tau$-invariant vector fields on $\widehat\Sigma$.

There is a differential operator denoted $\widehat D$ taking pairs of functions on $\widehat\Sigma$ to vector fields on $\widehat\Sigma$ lifting the operator $D$ acting on the analogous objects on $\Sigma$. We will describe this operator in more detail later in this section. One can check that $\widehat D\tau^*=\tau^*\widehat D$.
It follows that $\tau$ acts on the kernel and cokernel of $\widehat D$. Since
$\tau^2=1$ the only possible eigenvalues are $\pm 1$. Any vector field can be written as
$V=\frac12 (V+\tau\cdot V)+\frac12 (V-\tau\cdot V)$ where the first term is in the $+1$-eigenspace and the second is in the $-1$-eigenspace. Furthermore, there is a duality that carries one eigenspace to the other. We will further explain this duality later in this section. A similar thing holds for the pairs of functions. The important consequence to notice now is that the $+1$-eigenspace of the kernel of $\widehat\Sigma$ may be identified with $H^0(\Sigma)\oplus H^2(\Sigma)$ and this
is isomorphic to the $-1$-eigenspace. Similarly, the $+1$-eigenspace of the cokernel
of $\widehat\Sigma$ may be identified with $H^1(\Sigma)$ and this
is isomorphic to the $-1$-eigenspace. It follows from this and equation (\ref{indd}) that
\begin{equation}\label{inddh}
\begin{array}{l}\text{dim}(\text{ker}(\widehat D))=2\\
\text{dim}(\text{coker}(\widehat D))=2g\,.
\end{array}
\end{equation}

\begin{defn} The {\it index} of an operator $D$ is
\[
\text{\rm index}(D)=\text{\rm dim}(\text{\rm ker}(D))-\text{\rm dim}(\text{\rm coker}(D))\,.
\]
\end{defn}
Thus $\text{index}(\widehat D)=2-2g$. By counting the vertices, edges and faces of any polygonal decomposition of a surface of genus $g$, one can compute that such a surface has Euler characteristic equal to $2-2g$. Thus the index is the Euler characteristic and the Gauss-Bonnet theorem (equation \eqref{gb2}) gives:
\[
\text{index}(\widehat D)=\frac{1}{2\pi} \int_{{\widehat\Sigma}} K\,d\,\text{Area}\,.
\]

\begin{quote}
The structure of this formula generalizes. The index of any elliptic differential operator may be expressed as an integral over the base manifold.
\end{quote}

\section*{Rotation}
Before stating the index theorem, we will discuss rotation in a bit more detail.
To describe a family of vectors rotating along a path, let $x=x(t)$, $y=y(t)$ be parametric equations of the path, and let $V(t)=(u(t),v(t))$ be the vector at each point of the path. Let $V_0=(u_0,v_0)=(u(0),v(0))$ be the initial
vector in the family, and let $\theta=\theta(x(t),y(t))$ be the amount that the vector rotates from the initial position. In this case we may write,
\[
V=\begin{bmatrix}u \\ v\end{bmatrix}=\begin{bmatrix}\cos(\theta) & -\sin(\theta) \\ \sin(\theta) & \cos(\theta)\end{bmatrix}\begin{bmatrix}u_0 \\ v_0\end{bmatrix}\,.
\]
Taking the differential of this equation gives:
\[
\begin{aligned}
dV=\begin{bmatrix}du \\ dv\end{bmatrix} &=\begin{bmatrix}-\sin(\theta)(\theta_x dx+\theta_y dy) & -\cos(\theta)(\theta_x dx+\theta_y dy) \\ \cos(\theta)(\theta_x dx+\theta_y dy) & -\sin(\theta)(\theta_x dx+\theta_y dy)\end{bmatrix}\begin{bmatrix}u_0 \\ v_0\end{bmatrix} \\ \\ &=\begin{bmatrix}0 & -(\theta_x dx+\theta_y dy) \\ (\theta_x dx+\theta_y dy) & 0\end{bmatrix}\begin{bmatrix}\cos(\theta) & -\sin(\theta) \\ \sin(\theta) & \cos(\theta)\end{bmatrix}\begin{bmatrix}u_0 \\ v_0\end{bmatrix}\,.
\end{aligned}
\]
Thus we may write:
\[
dV+\omega V=0\,,
\]
where
\[
\omega = \begin{bmatrix}0 & (\theta_x dx+\theta_y dy) \\ -(\theta_x dx+\theta_y dy) & 0\end{bmatrix}\,.
\]

In the above description, it appears that the angle that the vector rotates depends only on the point where the vector is based, thus there would be a globally defined
rotation angle that would be a function of position. This need not be the case. The amount that the vector rotates might depend on the path that the vector follows. A good example to think about here is parallel transporting a vector around a $\pi/2 - \pi/2 - \pi/2$ triangle on a sphere. To model this situation one just replaces
$(\theta_x dx+\theta_y dy)$ by a differential $1$-form. A matrix of $1$-forms
is called a {\it connection form}. 

Now that we have this method to describe the rotation of a vector along a path, we are very close to a proof of the Gauss-Bonnet theorem. We just need to know a bit more about differential forms. We will take a short aside to describe the minimal amount of information we need about differential forms and to describe the operator $\widehat D$ on a closed surface.

It is beyond the scope of an expository article to describe differential forms in general, but we will say that a $k$-form is exactly the object that can be integrated over an oriented $k$-dimensional object to get a number. Thus a $1$-form in two dimensional space has an expression $\alpha = f(x,y)\,dx+g(x,y)\,dy$ and we can make sense of $\int_\gamma \alpha$ for any smooth path $\gamma$. Similarly, an example of a $2$-form would be $h(x,y)\,dx\wedge dy$. Integrating this over a surface patch in the plane would just give the area of the patch.

The space of $k$-forms on a manifold $M$ is denoted by
$C^\infty(\wedge^kM)$. The exterior derivatives map $k$-forms to $(k+1)$-forms.
In coordinates the exterior derivative is given by
\[
d(f\,dx^{i_1}\wedge d x^{i_2}\wedge\cdots dx^{i_k})=(df)\wedge dx^{i_1}\wedge d x^{i_2}\wedge\cdots dx^{i_k}\,,
\]
where $df$ is the usual differential. Note that $\wedge$ is an antisymmetric product. The differential of a function is very closely related to the gradient of a function. Similarly, the exterior derivative of a $1$-form is very closely related to the curl of a vector field. Indeed, it is the exterior derivative that is used to extend the operator $D$ from a subset of the plane to an operator $\widehat D$ on a surface $\Sigma$.

To describe the adjoint of the exterior derivative on $1$-forms, we will use an important duality on forms. Any $1$-dimensional volume element in $\R^3$ corresponds to a ``perpendicular" $2$-dimensional volume element in $\R^3$. Thus we can associate $*dz=dx\wedge dy$ to $dz$. More generally on an oriented $n$-dimensional Riemannian manifold we have the {\it Hodge star operator}:
\[
*:C^\infty(\wedge^kM)\to C^\infty(\wedge^{n-k}M)\,.
\]

With this background, the operator $\widehat D:C^\infty(\wedge^0\widehat\Sigma)\oplus C^\infty(\wedge^0\widehat\Sigma)\to C^\infty(\wedge^1\widehat\Sigma)$ is defined by
$\widehat D(u,w):=du-*dw$.

Returning to the proof of the Gauss-Bonnet theorem, we will make a few definitions and observations. The angle from a parallel translate of the initial tangent vector (to a curve) to the current tangent vector is called the exterior angle. We could also measure the angle from the current tangent (or normal) vector to a parallel translate of the initial tangent (or normal) vector. This would be negative of the exterior angle. We will denote this negative exterior angle by $\theta_{NT}$.

In our musings on rotating vectors we did not specify which coordinate system we were using. Of course specific values would change in different coordinate systems. One natural coordinate system to use on the boundary of a patch of surface would have the outer unit normal vector as the first coordinate vector and the tangent vector as the second. Thus $\theta_{NT}$ is the amount of rotation of a parallel translate of the initial normal vector in the normal-tangent coordinate system.

The normal-tangent coordinate system does not extend across the entire patch of a simply-connected surface, thus to relate rotation around the boundary of a surface patch to something inside the surface we should use a globally defined coordinate system on the patch. Given a global {\bf ij} coordinate system on the surface patch, we can relate the rotation angle of the parallel translate of the tangent vector in this coordinate system (denoted by $\theta_{\bf ij}$) to the angle in the normal-tangent coordinates by
\[
\theta_{\bf ij}=\theta_{NT}+\varphi\,.
\]

\begin{figure}[!ht]  
\hskip140bp
\includegraphics[width=55mm]{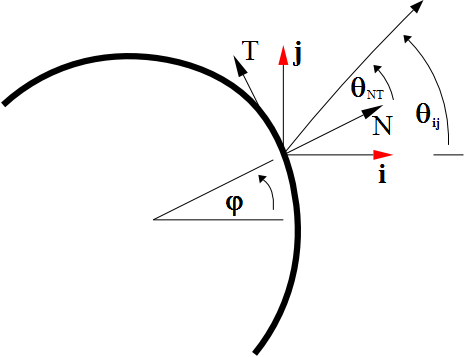}\caption{Angles in different coordinate systems}\label{coord}
\end{figure}

Our account of connection forms has been a bit like a description of linear algebra starting with matrices. A matrix is an expression for a linear map with respect to coordinate systems in the domain and codomain. Changing the coordinate systems changes the matrix, and there is a coordinate independent definition of a linear map. Similarly, there is a coordinate independent definition of a connection, and one can understand how the connection form expressions for a connection change when coordinates are changed. We do not need to go into any of these details here. However, more detail about all of these concepts may be found in the book by Singer and Thorpe \cite{ST}.

One traditional way to keep track of the exterior angle is via its derivative with respect to arclength $s$.
\begin{quote}
The geodesic curvature is the derivative of the exterior angle with respect to arclength:
\[
k_g=\frac{d\text{\rm ext ang}}{ds}\,.
\]
\end{quote}

In the case of two geodesics meeting at a corner, the geodesic curvature would be a singular measure. In any case the integral of the geodesic curvature over a curve is the exterior angle extended through the curve. Thus to compute the exterior angle around the boundary of patch $D$ in a surface we can write
\[
\begin{aligned}
2\pi - \text{\rm ext ang} &= 2\pi - \int_{\partial D} k_g\, ds \\
&= 2\pi + \int_{\partial D} d\theta_{NT} \\
&=2\pi + \int_{\partial D} d\theta_{\bf ij}-d\varphi \\
&=\int_{\partial D} d\theta_{\bf ij} \\
&=\int_{\partial D} \omega^1_2\,.
\end{aligned}
\]

Here $\omega^1_2$ is the form in the first row and second column of the connection form.
We will now use Stokes' theorem to write this as an integral over $D$. The result will be the integral of the sectional curvature over $D$. We pause to define the curvature and the sectional curvature.
\begin{quote}
The curvature is a matrix of $2$-forms obtained as a suitable derivative of the connection form.
\end{quote}
Specifically,
\[
F=d\omega+\omega\wedge\omega\,.
\]
Here the exterior derivative acts entry-wise on the connection form, and the wedge product is matrix multiplication with wedge products on the entries. When $\omega$ is an antisymmetric $2\times 2$ matrix, $\omega\wedge\omega=0$. In the case of an oriented surface, every $2$-form will be a multiple of the area form. Thus we can write $F^1_2=K\,d\text{Area}$, where $K$ is some function. This function is known as the {\it sectional curvature}. Returning to the Gauss-Bonnet theorem, we have
\[
\begin{aligned}
2\pi - \text{\rm ext ang} &=\int_{\partial D} \omega^1_2 \\
&= \int_D d\omega^1_2 = \int_D F^1_2 = \int_D K\,d\text{Area}\,.
\end{aligned}
\]


\section*{Theorem}
A result similar to the Gauss-Bonnet theorem holds for any elliptic differential operator. We state the result here and will explain terms used in the result as we encounter them in the later section on applications of the theorem.
\begin{thm}[Atiyah-Singer] If $E$ and $F$ are complex vector bundles over a closed oriented manifold $X$ and $D: C^\infty(E)\to C^\infty(F)$ is an elliptic differential operator, then
\[
\text{\rm index}(D)=\int_X (-1)^{n(n+1)/2} \text{\rm ch}(E-F)\text{\rm td}(T_\C X)\left(e(TX)\right)^{-1}\,.
\]
\end{thm}

Experts will know that a slightly simpler formula may be written by replacing the integral over the base manifold by an integral over the co-tangent space of the base manifold. We prefer this form because it fits with the nice slogan that ``the index of an elliptic differential operator may be expressed as the integral over something on the base."

We now sketch an argument that the index of an elliptic differential operator can be expressed as the integral of some expression that only depends on the geometry of the functions in the domain and codomain of the operator. This argument can be understood using only concepts from linear algebra.

The first step is
a simple computation using the definition of the adjoint that shows $\text{ker}(D)=\text{ker}(D^*D)$ and $\text{ker}(D^*)=\text{ker}(DD^*)$ provided that the inner products are positive definite. Notice that $D^*Dx=\lambda x$ implies that
$DD^*(Dx)=\lambda DX$, thus
\[
D:\text{ker}(D^*D-\lambda)\to\text{ker}(DD^*-\lambda)\,.
\]
If $\lambda\ne 0$ this is an isomorphism. Indeed, the same argument shows that $D^*$ takes $\text{ker}(DD^*-\lambda)$ to $\text{ker}(D^*D-\lambda)$ and
$D^*Dx=\lambda x$ implies that $\lambda^{-1}D^*$ is the inverse to $D$ on this subspace.

It follows that
\begin{equation}\label{supertrace}
\begin{array}{rl}
\text{index}(D)&=\text{dim}(\text{ker}(D))-\text{dim}(\text{ker}(D^*)) \\
&=\text{dim}(\text{ker}(D^*D))-\text{dim}(\text{ker}(D^*D)) \\
&=\sum_{\{u_\lambda\}}e^{-t\lambda}-\sum_{\{v_\mu\}}e^{-t\mu}\\
&=\text{Tr}(e^{-tD^*D})-\text{Tr}(e^{-tDD^*})\,.
\end{array}
\end{equation}
Here $\{u_\lambda\}$ and $\{v_\mu\}$ are orthonormal bases of eigenfunctions of $D^*D$ and $DD^*$ respectively. A straight forward exercise shows that distinct eigenspaces of $D^*D$ are perpendicular. It takes a bit more work to show that the eigenspaces are finite dimensional when $D$ is what is known as an elliptic operator and to show that the sums in the second to last line above converge. Notice that all of the terms in the sums corresponding to non-zero eigenvalues cancel because the $\lambda$ eigenspaces of $D^*D$ and $DD^*$ are isomorphic for all non-zero $\lambda$.
Notice this expression is independent of $t$, so we could use any $t$ value or take the limit as $t\to 0^+$.

The operator takes vector-valued functions to vector-valued functions. Given an inner product on the vector space one defines an inner product on the space of functions via
\[
\langle u, w\rangle_{C^\infty(E)}:=\int_X \langle u(x), w(x)\rangle_E\,d\,\text{vol}_X\,.
\]
A vector $u(y)$ can be considered as an operator taking $\R$ to $E$. The adjoint of this operator is denoted by $u(y)^\dagger$.

Setting $u_0(x)=\sum_{\{u_\mu\}}a_\mu u_\mu(x)$ one has
\begin{align*}
u(x,t)&:= e^{-tD^*D}u_0(x)\\
&=e^{-tD^*D}\sum_{\{u_\mu\}}a_\mu u_\mu(x)\\
&=\sum_{\{u_\mu\}}a_\mu e^{-\mu t} u_\mu(x)\\
&=\sum_{\{u_\mu\}}\sum_{\{u_\lambda\}}a_\mu e^{-\lambda t} u_\lambda(x)\int_Xu_\lambda(y)^\dagger u_\mu(y)\,d\,\text{vol}_y\\
&=\int_X \left(\sum_{\{u_\lambda\}}e^{-t\lambda}u_\lambda(x)u_\lambda(y)^\dagger\right)\sum_{\{u_\mu\}}a_\mu u_\mu(y)\,d\,\text{vol}_y\\
&=\int_Xk^{D^*D}_t(x,y)u_0(y)\,d\,\text{vol}_y\,,
\end{align*}
since the integral in the fourth line is equal to one when $\lambda=\mu$ and zero otherwise (as $\{u_\lambda\}$ is a orthonormal basis).
The expression
\[
k^{D^*D}_t(x,y):=\sum_{\{u_\lambda\}}e^{-t\lambda}u_\lambda(x)u_\lambda(y)^\dagger\,,
\]
is called the heat kernel of the operator  $D^*D$. It is apparent  from the first line that $\frac{\partial}{\partial t}u(x,t)=-D^*Du(x,t)$. Also $u(x,0)=u_0(x)$. In other words we see that $u(x,t)$ is the solution to this partial differential equation with initial data $u_0(x)$. This equation is called a heat equation because the scaler version of it models heat flow.

The $t\to 0^+$ limit of $k^{D^*D}_t(x,y)$ must be a delta function. In particular this only depends on the geometry of the space of functions and not on the specific form of the differential operator. Now,
\begin{align*}
\text{Tr}(e^{-tD^*D})&=\sum_{\{u_\lambda\}}e^{-t\lambda}\\
&=\sum_{\{u_\lambda\}}e^{-t\lambda}\int_Xu_\lambda(y)^\dagger u_\lambda(y)\,d\,\text{vol}_y\\
&=\sum_{\{u_\lambda\}}e^{-t\lambda}\int_X\text{Tr}(u_\lambda(y)u_\lambda(y)^\dagger) \,d\,\text{vol}_y\\
&=\int_X\text{Tr}(\sum_{\{u_\lambda\}}e^{-t\lambda}u_\lambda(y)u_\lambda(y)^\dagger) \,d\,\text{vol}_y\\
&=\int_X\text{Tr}(k^{D^*D}_t(x,y)) \,d\,\text{vol}_y\,.
\end{align*}

Combining this and the analogous expression for $DD^*$ with the expression for the index in equation (\ref{supertrace}) gives
\[
\text{index}(D)=\int_X\lim_{t\to 0^+}\left(\text{Tr}(k^{D^*D}_t(x,y))-\text{Tr}(k^{DD^*}_t(x,y))\right) \,d\,\text{vol}_y\,.
\]

To summarize, in this section we stated the amazing theorem of Atiyah and Singer relating the index of a differential operator to an integral over an expression depending only on the geometry of the space of functions in the domain and codomain of the operator. The outline that we gave just used elementary linear algebra. To identify the integrand and turn the outline into a rigorous proof there are a number of technical convergence questions that one must answer. The first outline of a proof of the index theorem by Atiyah and Singer used something called cobordism theory to reduce the proof to a check of several special cases \cite{AS1}. Their paper with a full proof, used a different technique known as K-theory \cite{AS2}. Their original cobordism proof was presented by Palais \cite{Pal}. A heat equation proof of the index theorem was conjectured by Kotake \cite{ko}, and proofs and refinements of the heat equation proof were given in a number of papers, \cite{G,ABP}. A simplified argument to identify the integrand was based on ideas about supersymmetry in physics leads to the short proof of the index theorem given by Getzler \cite{Ge}. There is now a vast literature on the index theorem and its generalizations with many excellent articles and books. Analytical details following the outline that we have given here appear for example in the textbook \cite{taylor}. It is time to look at some of the amazing and varied applications of the index theorem.

\section*{Jackpot}
You just won at a slot machine! The casino always empties the slot when it reaches \$$1,000.00$. The machine only has quarters and nickels. One possible jackpot is $2$ quarters and $3$ nickels, one possibility is $4,000$ quarters and there are many more possibilities. All possibilities totaling \$$1.00$ or less are displayed in the figure below.
\begin{figure}[!ht]  
\hskip70bp
\includegraphics[width=85mm]{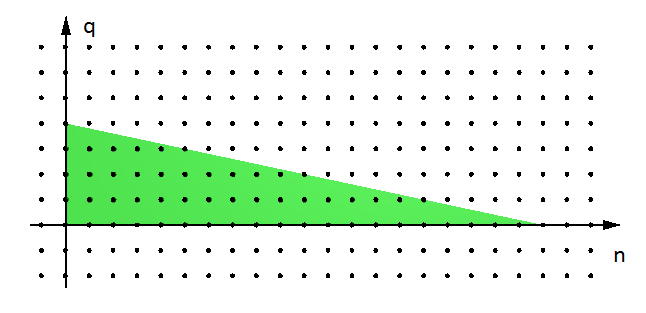}\caption{Jackpots no more than \$$1.00$}\label{jackpot}
\end{figure}
How many different jackpots are there?
The index theorem can answer this question. To see how, let $q$ denote the number of quarters and $n$ denote the number of nickels. The possibilities are simply the integer lattice points satisfying
\[
q\ge 0, \quad n\ge 0, \quad 0.25q+0.05n\le 1,000\,.
\]
The last equation can be rewritten as $c:=20,000-5q-n\ge 0$, so the number is the same as the number of non-negative integer lattice points $(q,n,c)$ satisfying
\[
5q+n+c=20,000=5k\,.
\]
Here $k=4,000$ represents the most quarters the slot machine will award. Of course we can generalize the problem and consider different maximums.

In order to apply the index theorem, notice that each triple of lattice points can be encoded as the coefficients of a monomial $\theta^q\eta^n\gamma^c$. These monomials will be a basis of the kernel of a differential operator that we will construct.

As a warm up consider the space of degree $k$ polynomials defined over the complex numbers. We can identify this space of polynomials as the kernel of a differential operator.
Consider two examples: $f(z)=z\bar z$ is not a polynomial, and $g(z)=z^2$ clearly is a polynomial. If we set $\bar\partial:=d\bar z\,\frac{\partial}{\partial \bar z}$ we will obtain
$\bar\partial f = z\,d\bar z$ and $\bar\partial g=0$. A standard result from complex analysis implies that any function $h(z)$ satisfying $\bar\partial h=0$ may be expressed as a convergent power series $h(z)=\sum_{n=0}^\infty c_nz^n$.

To restrict to degree $k$ functions we need to further specify the geometry of the domain of the operator. First recall that the
version of the index theorem that we stated is valid for {\it closed} manifolds. To use the index theorem in this setting we add a point at infinity viewed as the ratio $[1:0]$.
The space of all ratios is
complex projective space
\[
\CP^1:=\{(z,w)\in \C^2-\{(0,0)\}\}/\sim\,,\quad (z,w)\sim (\lambda z,\lambda w)\,,
\]
for $\lambda\in \C-\{0\}$,  and we denote the equivalence class of $(z,w)$ (the ratio) by $[z:w]$.

The degree $k$ condition enters through the following definition:
\[
L_k:=\{(z,w,f)\in (\C^2-\{(0,0)\})\times\C\}/\sim\,,\quad (z,w,f)\sim (\lambda z,\lambda w,\lambda^k f)\,,
\]
for $\lambda\in \C-\{0\}$. Notice that we have a projection $\pi:L\to \CP^1$. The inverse image of each point has the structure of a vector space. This is an example of something called a {\it vector bundle}. The domain of the $\bar\partial$ operator, $C^\infty(L)$, is the set of smooth functions $s:\CP^1\to L$ such that $\pi\circ s=\text{id}_{\CP^1}$. Such sections may be written $s([z:w])=[z,w,f(z,w)]$ where $f(\lambda z,\lambda w)=\lambda^kf(z,w)$ and are completely specified by $s([z:1])=[z,1,f(z,1)]$ and $s([1:w])=[1,w,f(1,w)]$. The  operator is defined by $\bar\partial s([z:1])=[z,1, \frac{\partial f}{\partial \bar z}|_{(z,1)}d\bar z]$ and  $\bar\partial s([1:w])=[1,w, \frac{\partial f}{\partial \bar w}|_{(1,w)}d\bar w]$. Since $f(z,1)=z^kf(1,z^{-1})$ we can set  $z=w^{-1}$ for $z\ne 0$ and compute
\[
\begin{aligned}
\bar\partial s([z:1])&=[z,1, \frac{\partial f}{\partial \bar z}|_{(z,1)}d\bar z]\\
&=
[z,1, z^k\frac{d\bar w}{d\bar z}\frac{\partial f}{\partial \bar w}|_{(1,w)}d\bar z]\\
&=
[1,z^{-1}, \frac{d\bar w}{d\bar z}\frac{\partial f}{\partial \bar w}|_{(1,w)}d\bar z]\\
&=[1,w, \frac{\partial f}{\partial \bar w}|_{(1,w)}d\bar w]=\bar\partial s([1:w])\,,
\end{aligned}
\]
to see that the operator is well defined. For $s\in\text{ker}(\partial)$ we have
$\frac{\partial f}{\partial \bar w}|_{(1,w)}=0$, so $f(1,w)=\sum_{n=0}^\infty b_nw^n$
and $f(z,1)=z^kf(1,z^{-1})=z^k \sum_{n=0}^\infty b_nz^{-n}$. The only way this latter expression can can be a power series is if $b_n=0$ for $n>k$ in which case it is a polynomial of degree $k$.

The collection of all differentials forms a vector bundle known as the cotangent bundle
$T^*\CP^1$ and the codomain of the $\bar\partial$ operator is a combination of the cotangent bundle and the line bundle known as the tensor product $T^*\CP^1\otimes L$. One can show that the cokernel of $\bar\partial$ is zero so that its index is just the dimension of the space of degree $k$ polynomials. It follows that
\[
\text{index}(\bar\partial:C^\infty(L_k)\to C^\infty(T^*\CP^1\otimes L_k)) = k+1\,.
\]

We can apply the index theorem to compute this. In the process, we will learn the definitions of all of the terms in the index theorem. The first thing to notice is that all of the vector bundles in the index theorem are complex vector bundles. We can use the complex numbers to model a rotating vector in the plane similar to the way that we used vectors to model rotation in the plane in the rotation section. Thus we can set $V=e^{{\bf i}\theta}V_0$ and repeat the motivation from the first section to see that rotation along paths may be modeled by a purely imaginary $1$-form, $\omega$. For higher dimensional complex vector spaces rotations can be modeled by skew-hermitian matrices with $1$-form entries.

It this toy example, we can take  a different model for $L_k$:
\[
L_k=\{(z,w,f)| |z|^2+|w|^2 = k \}/\sim\,,
\]
where $(z,w,f)\sim (\lambda z, \lambda w \lambda^k f)$ for $\lambda\in S^1$. In this model, it is easy to see that
the following is a connection form on $L_k$:
\[
\omega_k=\text{Im}\left(-k\frac{\bar z \,dz +\bar w \,dw}{|z|^2+|w|^2}  \right)\,.
\]
Since $\omega_k\wedge\omega_k=0$ in this case, the curvature will just be $F_k=d\omega_k$. Notice that $F_k=kF_1$.

To apply the index theorem to the toy problem and the jackpot problem, we will first use a simplification of the integrand that takes place over complex manifolds, in particular for what are known as twisted $\bar\partial$ operators one has:
\begin{equation}\label{RR}
(-1)^{n(n+1)/2} \text{\rm ch}(E-F)\text{\rm td}(T_\C X)(\left(e(TX)\right)^{-1}=\text{ch}(L)\text{td}(TX)\,.
\end{equation}
The heuristics here are that the Todd class is defined for complex vector bundles. When the base $X$ is a complex manifold, the tangent space will already be complex, so one can relate the Todd class of the complexification of the tangent bundle to the Todd class of the tangent bundle. This is explained nicely in the book by Shanahan \cite{S}.

By definition, the Chern character of a complex vector bundle $E$ is given by \[\text{ch}(E)=\text{tr}(\text{exp}(-\frac{1}{2\pi i}F))\,,\] where $F$ is the curvature of a connection form on the vector bundle. In the particular case of $L_k$ we can expand the exponential in a Taylor expansion to get
\[\text{ch}(L_k)=1-\frac{1}{2\pi i}F_k\,.\] Recall that the curvature is a matrix of $2$-forms, so a term with $F\wedge F$ would be a $4$-form, but $\CP^1$ is a real $2$ dimensional surface, so $F^{\wedge 2}$ and all higher powers of $F$ vanish on it.

By definition, the Todd class of a complex vector bundle $E$ is given by
\[
\text{td}(E)=\text{det}(-\frac{1}{2\pi i}F(1-\text{exp}(\frac{1}{2\pi i}F))^{-1})\,,
\]
where $F$ is the curvature of a connection on $E$. It turns out that the tangent bundle $T\CP^1$ (i.e. space of all tangent vectors to $\CP^1$) is isomorphic to the bundle $L_{2}$. To see this draw a vector field on $\CP^1$ viewed as a $2$-sphere, such a vector field may be drawn with exactly two non-degenerate zeros, and this corresponds to the fact that a generic degree $2$ polynomial has exactly two non-degenerate zeros. Expanding the expression for the Todd class in a Taylor expansion gives:
\[
\text{td}(T\CP^1)=1-\frac12\frac{1}{2\pi i}F_2.
\]

It follows that
\[
\text{ch}(L_k)\text{td}(T\CP^1) = 1-\frac{1}{2\pi i}F_k-\frac12\frac{1}{2\pi i}F_2
=1-(k+1)\frac{1}{2\pi i}F_1\,.
\]
Here we use the fact that the higher order terms vanish since this lives on a real $2$-dimensional surface, and $F_k=kF_1$.

We can complete the computation by evaluating $-\frac{1}{2\pi i}\int_{\CP^1} F_1$. The first thing to notice is that removing the point $[1:0]$ does not effect the value of the integral. Thus we can just work in the $w=1$ chart:
\[
\begin{aligned}
-\frac{1}{2\pi i}\int_{\CP^1} F_1 &=-\frac{1}{2\pi i}\int_{\C} d\omega_1\\
 &=-\lim_{R\to\infty}\frac{1}{2\pi i}\int_{|z|\le R} d\omega_1\\
 &=-\lim_{R\to\infty}\frac{1}{2\pi i}\int_{|z|= R} \omega_1\\
  &=\lim_{R\to\infty}\frac{1}{2\pi i}\int_{|z|= R} \text{Im}\left(\frac{\bar z \,dz}{|z|^2+1}  \right)\\
    &=\lim_{R\to\infty}\frac{1}{2\pi i}\frac{R^2}{R^2+1}\int_{|z|= R} \text{Im}\left(\frac{dz}{z}  \right)=1\,.
\end{aligned}
\]

Now return to the jackpot problem (counting the non-negative integers solutions to $5q+n+c=20,000=5k$) and consider the space
\[
\begin{array}{c}
L:=\{(\theta,\eta,\gamma,f) | 5|\theta|^2+|\eta|^2+|\gamma|^2=5k\}/\sim\,,\\ (\theta,\eta,\gamma,f)\sim (\lambda^5\theta,\lambda\eta,\lambda\gamma,\lambda^{5k}f)\,,
\end{array}
\]
for $|\lambda|=1$. There is a natural projection from $L$ to
\[
\begin{array}{c}
X:=\{(\theta,\eta,\gamma) | 5|\theta|^2+|\eta|^2+|\gamma|^2=5k\}/\sim\,,\\ (\theta,\eta,\gamma)\sim (\lambda^5\theta,\lambda\eta,\lambda\gamma)\,,
\end{array}
\]
for $|\lambda|=1$. Notice that the values of $(|\theta|^2,|\eta|^2,|\gamma|^2)$ recover the original polytope. For $(\theta,\eta,\gamma,f)\in\C^4$ we obtain an interesting space from which we can recover the lattice points in an interesting way. In this case we have,
\[
\begin{array}{c}
L\cong \{(\theta,\eta,\gamma,f) \in \C^4-\{0\}\}/\sim\,,\\ (\theta,\eta,\gamma,f)\sim (\lambda^5\theta,\lambda\eta,\lambda\gamma,\lambda^{5k}f)\,,\quad\text{\rm for} \ \lambda\in\C-\{0\}\,.
\end{array}
\]
A section in $C^\infty(L)$ is expressed as $s([\theta,\eta,\gamma])=[\theta,\eta,\gamma,f(\theta,\eta,\gamma)]$ for some smooth $f$ satisfying $f(\lambda^5\theta,\lambda\eta,\lambda\gamma)=\lambda^{5k}f(\theta,\eta,\gamma)$.


As in the toy problem of polynomials of degree $k$, good examples of such smooth functions are given by polymomials. The monomial $f(\theta,\eta,\gamma)=\theta^6\eta^2\gamma^3$ satisfies $f(\lambda^5\theta,\lambda\eta,\lambda\gamma)=\lambda^{5\cdot 7}f(\theta,\eta,\gamma)$.
We can think of it as a representation of a pay out of $6$ quarters, $2$ nickels, on a slot machine with a maximum pay out of $7$ quarters. Thus it is $3$ nickels short of the maximum pay out. Compare these numbers with the exponents of the specific monomial $f$. It should be clear that we can create a bijection between monic monomials and lattice points in the lattice polytope describing all possible jackpots.

In fact, a bit of thought should convince the reader that there is a bijection from the lattice points in a lattice polytope and the lattice points of the intersection of an affine subspace with the non-negative cone in a higher dimensional lattice. Just arrange each inequality defining the polytope to be of the form $\text{\it something}\ge 0$ and then declare each {\it something} to be a new (slack) variable. The lattice points in this latter section can then be identified with monic weighted-homogenous polynomials. If these polynomials can be identified as a basis of the kernel of an elliptic differential operator, then the Atiyah-Singer Theorem can be used to compute the number of lattice points in the lattice polytope. This is indeed the case. The reason for the construction of the complicated spaces $L$ and $X$ is so we can define these operators. We now construct the appropriate differential operator in the example of the jackpot problem.

Defining an operator on the sections $C^\infty(L)$ (in the $\gamma=1$ case) by
\[\bar\partial s = [\theta,\eta,1,\frac{\partial}{\partial\bar\theta}fd\bar\theta+\frac{\partial}{\partial\bar\eta}fd\bar\eta]\,,\]
with similar expressions for the $\theta=1$ and $\eta=1$ cases produces a well-defined operator in much the same way as the $\bar\partial$ operator from our toy example. Functions in the kernel of this operator will be holomorphic, and will thus have power series representations. In order for these power series to match up in any chart ($\gamma=1$, $\eta=1$, or $\theta=1$) the corresponding section will have to be a suitably weighted-homogeneous polynomial. Thus the dimension of this kernel is exactly the number of lattice points in the lattice polygon. It takes a bit of work to show that the cokernel of this operator is trivial. Once this is done we see that the number of lattice points is just $\text{index}(\bar\partial)$.

In the jackpot example, the integral arising in the index theorem is fairly complicated.
Luckily, the spaces in question are very symmetric and this symmetry may be used to simplify the integral. Let $T^2:=\{(\lambda_1,\lambda_2)\in\mathbb{C}^2 | |\lambda_1|=1, |\lambda_2|=1\}$ be the two dimensional torus. It acts on $L$ as follows:
\[
(\lambda_1,\lambda_2)\cdot [\theta,\eta,\gamma,f] := [\lambda_1\theta,\lambda_2\eta,\gamma,f]\,.
\]
There is a similar action on $X$. The relevant (localization) theorem is that for any equivariant class $\phi$ one has:
\begin{equation}\label{localize}
\int_M\phi = \sum_F\int_F
\frac{\iota_F^*\hat\phi}{e(N(F))}\,.
\end{equation}
Here the sum is taken over all components of the fixed point set, $F$.
Justifying this nice formula of Atiyah and Bott \cite{AB} would take us a bit too far afield, but we will explain how to use it.

The first thing to do is to compute the fixed point set of the $T^2$ action on $X$. For
$[\theta,\eta,\gamma]$ to be fixed by all elements of $T^2$, there needs to be a $\lambda$ such that $(\lambda_1\theta,\lambda_2\eta,\gamma) = (\lambda^5\theta,\lambda\eta,\lambda\gamma)$. A bit of thought shows that the only possibilities are $[1,0,0]$, $[0,1,0]$ or $[0,0,1]$. Thus the localization formula (\ref{localize}) will reduce the complicated integral in this application to an integral over a these three points, i.e. a sum of three terms. These terms must keep track of the
group action. The way this is accomplished is via the corresponding representations.

Every representation of $T^2$ decomposes into a direct sum of $1$-dimensional representations. Every $1$-dimensional representation takes the form $(\lambda_1,\lambda_2)\cdot z = \lambda_1^{-n}\lambda_2^{-m}z$ for some integers $n$ and $m$. We denote these representations by $n\alpha_1+m\alpha_2$. The reason for the signs here comes from the proof of localization theorem \cite{AB}. We can now compute the induced $T^2$ representations on the tangent space of $X$ at the three fixed points as well as the induced representations on the fibers of $L$ over the three fixed points.

Using the complex description of $X$, we find a large open set around $[0,0,1]$ given by all points of the form $[\theta,\eta,1]$ with $(\theta,\eta)\in \mathbb{C}^2$. This open set is isomorphic to the tangent space, and we see that the corresponding representation is just $T_{[0,0,1]}X=(-\alpha_1)\oplus(-\alpha_2)$. Similarly, the action on the fiber of $L$ is trivial so $L_{[0,0,1]}=0$. Around $[0,1,0]$ we have a large open set consisting of all points of the form $[\theta,1,\gamma]$. The action is given by
\[
(\lambda_1,\lambda_2)\cdot [\theta,1,\gamma] = [\lambda_1\theta,\lambda_2,\gamma]
= [\lambda_2^{-5}\lambda_1\theta,1,\lambda_2^{-1}\gamma]\,.
\]
Here we used the definition of the equivalence relation to simplify the form of the representative, so we could identify the representation as $T_{[0,1,0]}X=(5\alpha_2-\alpha_1)\oplus(\alpha_2)$. A similar computation for $L$ shows that as a representation the fiber is $L_{[0,1,0]}=5k\alpha_2$.

Repeating the same argument for the point $[1,0,0]$ brings us to some technicalities. This time we get an open set around the fixed point modeled on $[1,\eta,\gamma]$ where
$(\eta,\gamma)\in\mathbb{C}^2$ however some of these points are equivalent. In particular, we have $[1,\eta,\gamma] = [1, \zeta\eta,\zeta\gamma]$ when $\zeta$ is any fifth root of unity. Thus the space $X$ has a singularity at this point and is not actually a manifold. The singularity here is not too bad. It is just the quotient of a Euclidean space by a finite group action (in this case $\mathbb{Z}_5$). An orbifold is a space that is build from patches that are quotients of Euclidean space by finite group actions. There are many generalizations of the Atiyah-Singer index theorem, in particular there is one to the orbifold situation.

\begin{quote}To integrate a measure over an orbifold one considers various coordinate patches, and lifts the corresponding measures to invariant ones on the corresponding Euclidean spaces and then divides by the order of the local symmetry group.
\end{quote}

We will follow this procedure here, we will also see that the representations used in the localization theorem need to be generalized a bit. Following the procedure that we used for the first two fixed points, we see that the action is given by
\[
(\lambda_1,\lambda_2)\cdot [1,\eta,\gamma] = [\lambda_1,\lambda_2\eta,\gamma]
= [1,\lambda_1^{-1/5}\lambda_2\eta,\lambda_1^{-1/5}\gamma]\,.
\]
Thus the tangent space at this point gives the representation
$T_{[1,0,0]}X=(\frac{1}{5}\alpha_1-\alpha_2)\oplus(\frac{1}{5}\alpha_1)$. Including the $f$ component in the above computation identifies the representation on the fiber as
$L_{[1,0,0]}=k\alpha_1$. Summarizing the computation so far we have:
\begin{equation}\label{Lrep}
\begin{aligned}
L_{[1,0,0]}&=k\alpha_1\,, \\
L_{[0,1,0]}&=5k\alpha_2\,, \\
L_{[0,0,1]}&=0\,,
\end{aligned}
\end{equation}
and
\begin{equation}\label{Trep}
\begin{aligned}
T_{[1,0,0]}X&=(\frac{1}{5}\alpha_1-\alpha_2)\oplus(\frac{1}{5}\alpha_1)\,, \\
T_{[0,1,0]}X&=(5\alpha_2-\alpha_1)\oplus(\alpha_2)\,, \\
T_{[0,0,1]}X&=(-\alpha_1)\oplus(-\alpha_2)\,.
\end{aligned}
\end{equation}

Since the fixed point set of the $T^2$ action is zero dimensional, the normal bundle to the fixed point set is just the tangent space to the base restricted to this set. We now need to compute the characteristic classes that appear in the localization formula. It is useful to introduce the Chern classes to do this. The total Chern class of a complex vector bundle $E$ is
\[
c(E):=\text{det}(\mathbb{I}-\frac{1}{2\pi i}F)\,,
\]
where $F$ is the curvature of any connection on $E$. All terms of this expression have even degree since $F$ is a matrix of $2$-forms. The part of the total Chern class with degree $2k$ is known as the $k$-th Chern class $c_k(E)$.

The top Chern class of a complex vector bundle is often defined to be the Euler class of the associated real bundle. We will turn this around. The Euler class $e(N(F))$ is the top
Chern class, and this is just the top degree symmetric polynomial i.e. the product of the eigenvalues of the normalized curvature matrix. These eigenvalues are nothing more than the irreducible factors of each representation. Thus,
\begin{equation}\label{en}
\begin{aligned}
e(N(F)_{[1,0,0]})&=(\frac{1}{5}\alpha_1-\alpha_2)(\frac{1}{5}\alpha_1)\,, \\
e(N(F)_{[0,1,0]})&=(5\alpha_2-\alpha_1)(\alpha_2)\,, \\
e(N(F)_{[0,0,1]})&=(-\alpha_1)(-\alpha_2)\,.
\end{aligned}
\end{equation}

Recall that the Chern character is given by $\text{ch}(E)=\text{tr}(\text{exp}(-\frac{1}{2\pi i}F))=\sum e^{x_i}$ where $x_i$ are the eigenvalues of the curvature matrix. Since the bundle $L$ is has one dimensional fibers, the normalized curvature matrix is a $1\times 1$ matrix and the eigenvalue of this matrix is just the representation! Thus (\ref{Lrep}) gives,
\begin{equation}\label{ch}
\begin{aligned}
\text{ch}(L_{[1,0,0]})&=1+k\alpha_1+\frac{1}{2}k^2\alpha_1^2\,, \\
\text{ch}(L_{[0,1,0]})&=1+5k\alpha_2+\frac{25}{2}k^2\alpha_2^2\,, \\
\text{ch}(L_{[0,0,1]})&=1\,.
\end{aligned}
\end{equation}

These are clearly polynomial in $k$, and as a result the entire integrand appearing in the index theorem will be a polynomial in $k$ as well. We will group the contributions to each term in the localization formula by the power of $k$.

In general the Todd class of a vector bundle is given by $\text{det}(-\frac{1}{2\pi i}F(1-\text{exp}(\frac{1}{2\pi i}F))^{-1})=\prod\frac{x_i}{1-e^{-x_i}}$ where $x_i$ are the eigenvalues of the normalized curvature matrix. For a rank two vector bundle over a four dimensional base this becomes
\[
\begin{aligned}
\frac{x_1}{1-e^{-x_1}}\frac{x_2}{1-e^{-x_2}}&=(1-x_1/2+x_1^2/6)^{-1}(1-x_2/2+x_2^2/6)^{-1} \\&=(1+x_1/2-x_1^2/12)(1+x_2/2-x_2^2/12)\\
&=1+\frac{1}{2}(x_1+x_2)+\frac{1}{12}(x_1^2+x_2^2+3x_1x_2)\\
&=1+\frac{1}{2}(x_1+x_2)+\frac{1}{12}(x_1+x_2)^2+\frac{1}{12}x_1x_2\\&=1+\frac{1}{2}c_1+\frac{1}{12}c_1^2+\frac{1}{12}c_2\,.
\end{aligned}
\]

The version of the index theorem for twisted $\bar\partial$ operators (\ref{RR}) and the Atiyah-Bott localization formula (\ref{localize}) tell us that
\[
\text{index}(\bar\partial)=\int_X\text{ch}(L)\text{td}(TX) = \sum_{p=[1,0,0],[0,1,0],[0,0,1]} \frac{\text{ch}(L_{p})\text{td}(T_{p}X)}{e(N(F)_{p})}\,.
\]
Certainly the highest degree terms of $\text{ch}(L_{p})\text{td}(T_{p}X)$ are $$\text{td}(T_{p}X)_4+\text{ch}(L_{p})_2\text{td}(T_{p}X)_2+\text{ch}(L_{p})_4$$
where we write the Chern character and Todd class as a sum of terms having different degrees, $\text{ch}(L_{p})=1+\text{ch}(L_{p})_2+\text{ch}(L_{p})_4$, and $\text{td}(T_{p}X)=1+\text{td}(T_{p}X)_2+\text{td}(T_{p}X)_4$.

To get the terms with no factor of $k$ in the localization formula we just need the highest degree terms of the Todd class. Combining the expressions for the Euler classes in (\ref{en}) with the expressions for the representations associated to the tangent bundle in (\ref{Trep}) and the highest degree term of the rank two Todd class gives
\[
\begin{aligned}
\frac{\text{td}(T_{[0,0,1]}X)_4}{e(N(F)_{[0,0,1]})}&=\frac{\left[\alpha_1^2+\alpha_2^2+3\alpha_1\alpha_2\right]}{12\alpha_1\alpha_2}\,,\\
\frac{\text{td}(T_{[0,1,0]}X)_4}{e(N(F)_{[0,1,0]})}&=\frac{\left[(5\alpha_2-\alpha_1)^2+\alpha_2^2+3(5\alpha_2-\alpha_1)\alpha_2\right]}{12(5\alpha_2-\alpha_1)\alpha_2}\,,\\
\frac{\text{td}(T_{[1,0,0]}X)_4}{e(N(F)_{[1,0,0]})}&=\frac{1}{5}\frac{\left[(\frac{1}{5}\alpha_1-\alpha_2)^2+(\frac{1}{5}\alpha_1)^2 +3(\frac{1}{5}\alpha_1-\alpha_2)(\frac{1}{5}\alpha_1)\right]}{12(\frac{1}{5}\alpha_1-\alpha_2)(\frac{1}{5}\alpha_1)}\,.
\end{aligned}
\]
The factor of $1/5$ in the last line is there because this term comes from the orbifold chart with non-trivial local symmetry group, and we have to divide by the order of this symmetry group. It is an exercise in algebra to show that the sum of these expressions is $1$.

To get the terms with a factor of $k$ in the localization formula we just combine the degree two portion of the Chern character with the degree two portion of the Todd class and the Euler class (\ref{ch}, \ref{Trep}, \ref{en}):
\[
\begin{aligned}
\frac{\text{ch}(L_{[0,0,1]})_2\text{td}(T_{[0,0,1]}X)_2}{e(N(F)_{[0,0,1]})}&=0\,,\\
\frac{\text{ch}(L_{[0,1,0]})_2\text{td}(T_{[0,1,0]}X)_2}{e(N(F)_{[0,1,0]})}&=\frac{5k\alpha_2\left[(5\alpha_2-\alpha_1)+\alpha_2\right]}{2(5\alpha_2-\alpha_1)\alpha_2}\,,\\
\frac{\text{ch}(L_{[1,0,0]})_2\text{td}(T_{[1,0,0]}X)_2}{e(N(F)_{[1,0,0]})}&=\frac{1}{5}\frac{k\alpha_1\left[ (\frac{1}{5}\alpha_1-\alpha_2)+(\frac{1}{5}\alpha_1) \right]}{2(\frac{1}{5}\alpha_1-\alpha_2)(\frac{1}{5}\alpha_1)}\,.
\end{aligned}
\]
Algebra shows that the sum of these terms is $\frac{7}{2}k$.

To get the terms with a factor of $k^2$ in the localization formula we just combine the degree four portion of the Chern character with  the Euler class (\ref{ch}, \ref{en}):
\[
\begin{aligned}
\frac{\text{ch}(L_{[0,0,1]})_4}{e(N(F)_{[0,0,1]})}&=0\,,\\
\frac{\text{ch}(L_{[0,1,0]})_4}{e(N(F)_{[0,1,0]})}&=\frac{25k^2\alpha_2^2}{2(5\alpha_2-\alpha_1)\alpha_2}\,,\\
\frac{\text{ch}(L_{[1,0,0]})_4}{e(N(F)_{[1,0,0]})}&=\frac{1}{5}\frac{k^2\alpha_1^2}{2(\frac{1}{5}\alpha_1-\alpha_2)(\frac{1}{5}\alpha_1)}\,.
\end{aligned}
\]
Algebra shows that the sum of these terms is $\frac{5}{2}k^2$. Thus the number of lattice points in the generalized lattice polygon for the jackpot problem is
\[
\frac{5}{2}k^2+\frac{7}{2}k+1\,.
\]
One can verify directly that this formula is correct for $k=0$, $1$, and $2$. Together with the expected quadratic growth in $k$ this gives an easier derivation of the specific formula. A different way to derive this formula would be to use Pick's theorem relating the area, and perimeter of a lattice polygon to the number of lattice points inside of it, \cite{pick}.

It is clear that the Atiyah-Singer index theorem is not the easiest way to solve this problem. The point is that the index theorem has a large number of surprising and varied applications. It may often be used to deduce unexpected connections between different types of mathematics, and often points to different generalizations of results. The discussion so far, should make it clear that there are higher dimensional versions of Pick's theorem that follow from the index theorem. Pommersheim wrote a nice dissertation in the mid 1990s about the application of the index theorem to lattice counting questions. Also see \cite{p3,p1,p2}.
We obtain the answer to our original jackpot question by substituting $k=4,000$ into the general formula to see that there are $40,014,001$ possible jackpots of \$$1,000$ or less with nickels and quarters. This is a large number, but it may be smaller than the number of possible applications of the index theorem.

\section*{Knotty}\label{knotty}
Things get a bit knotty now. The index theorem has been used to establish many results in knot theory and geometric topology. A {\it knot} is a smooth embedding of a circle into $\mathbb{R}^3$ (one glues the ends of the string together to keep the knot from coming undone). By scaling we may assume that the image of the knot is contained in $[-1,1]^3$. A knot is trivial if it is the boundary of a smooth disk in $\mathbb{R}^3$.

\begin{figure}[!ht]  
\hskip70bp
\includegraphics[width=85mm]{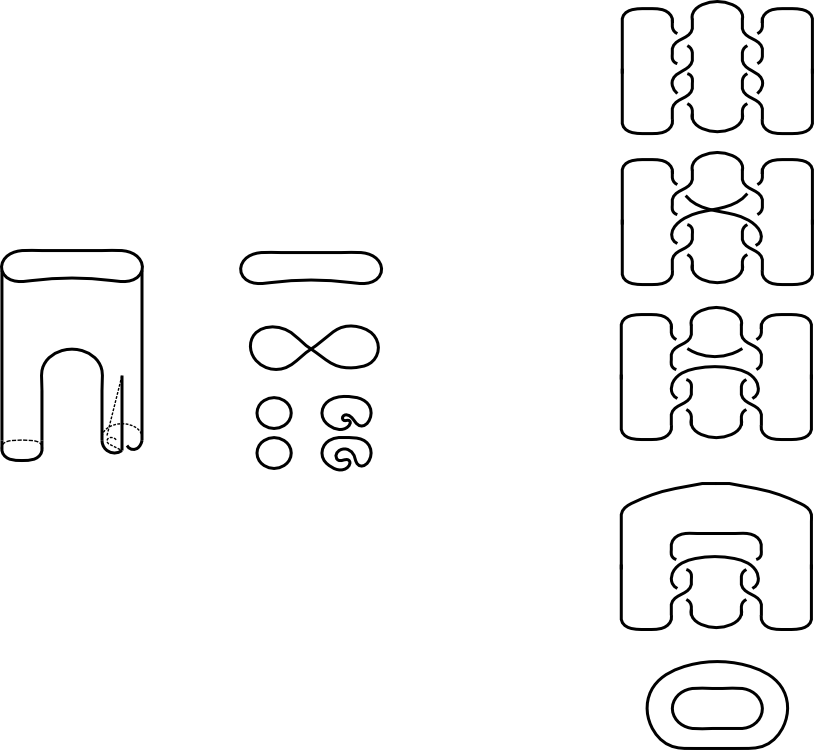}\caption{Movie of a slice disk}\label{gran}
\end{figure}

A non-trivial knot may still bound a smooth disk in $[-1,1]^4$.  Knots with this property are called slice. An example of a slice knot together with a movie of the
slice disk is displayed in figure \ref{gran}.
The left side has a movie representation of a pair of pants embedded in $[-1,1]^3$ to make it into a disk a "sock" can be sewed on the the bottom of each pant leg. The right side has the analogous movie picture in the four dimensional case. The top frame of this movie is the granny knot, and the entire movie displays the pants that granny wears.

\begin{figure}[!ht]  
\hskip150bp
\includegraphics[width=29mm]{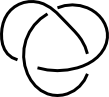}\caption{The trefoil knot}\label{t23}
\end{figure}

There are algebraic invariants that can demonstrate that certain knots (including the trefoil) are not slice, however,
we will use the index theorem to prove that the trefoil knot in figure \ref{t23} is not slice. There are easier ways to prove that this knot is not slice, but the ideas that we use are useful for more complicated knots and the index theorem is a very effective tool to study such problems. In particular one of the early applications of the index theorem to topology was the work by Casson and Gordon establishing that there are algebraically slice knots that are not slice, \cite{CG}. We are essentially using their argument here. More recently Jae Choon Cha has used these ideas to obtain very subtle invariants of concordance in the formidable paper \cite{cha}.

\begin{figure}  
\hskip55bp
\includegraphics[width=100mm]{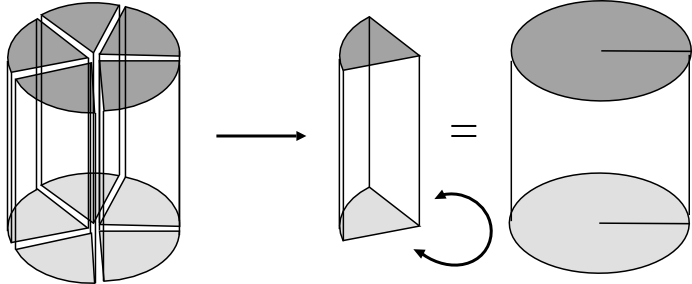}\caption{A branched cover}\label{branch}
\end{figure}

The proof that the trefoil is not slice is an indirect proof. Assuming that it is the boundary of a smooth disk, we construct a special type of manifold and apply the index theorem to get a contradiction. The construction begins with something called a  $5$-fold branched cover. A branched cover is a map from one manifold to another that is locally modeled on the map $\mathbb{R}^k\times\mathbb{C}\to\R^k\times\C$ given by
$(v,z)\mapsto (v,z^d)$. Figure \ref{branch} displays a $d=5$-fold branched cover of
$D^3\cong D^2\times[-1,1]$ to $D^3\cong D^2\times[-1,1]$ branched along an embedded $D^1$. This is a good schematic for the $5$-fold branched cover branched over the assumed slice disk. We will call this hypothetical cover $R$.

At the end of this section, we will show that the $5$-fold branched cover of the three sphere $S^3$ branched along the trefoil knot is the boundary of a spin manifold with signature $-8$, denoted $V_{\text{\rm res}}$. We will explain the signature and the notion of a spin manifold soon. The important points up to here are that the hypothetical cover $R$ has no $2$-cycles, so the manifold $X=\bar V_{\text{\rm res}}\cup R$ will be a closed, spin $4$-manifold with signature $-8$.

We now describe the signature of a $4$-dimensional manifold. In $4$-dimensions the Hodge star operator sends a $2$-form $dx^1\wedge dx^2$ to the ``perpendicular" form $dx^3\wedge dx^4$. The second cohomology,
\[
H^2(X) = \frac{\text{\rm Ker}(d:C^\infty(\wedge^2X)\to C^\infty(\wedge^3X))}{\text{\rm Im}(d:C^\infty(\wedge^1X)\to C^\infty(\wedge^2X))}\,,
\]
can then be decomposed into
positive and negative subspaces $H^2_+(X)$ and $H^2_-(X)$. By definition the signature of $W$ is $\text{Sign}(X):=\text{dim}(H^2_+(X))-\text{dim}(H^2_-(X))$.
This looks like the index of a differential operator. This is indeed the case, the signature is the index of the signature operator. An application of the index theorem gives
\[
\text{Sign}(X)=\frac{1}{3}\int_X p_1(X)\,.
\]

In four dimensions one may use the quaternions to describe spin structures. The quaternions are analogous to the complex numbers. They can be constructed by adjoining three elements, ${\bf i}$, ${\bf j}$, and ${\bf k}$ to the real numbers such that ${\bf i}^2={\bf j}^2={\bf k}^2={\bf ijk}=-1$.
A spin structure on a $4$-manifold consists of a pair of quaterionic vector bundles, $W_+$ and $W_-$ together with so called clifford multiplication operations $c:T^*X \times W_\pm \to W_\mp$. It is known that a simply-connected $4$-manifold admits a spin structure exactly when the self intersection number of every cycle is even. One defines an elliptic differential operator (the Dirac operator) $\partial\!\!\!/:C^\infty(W_+) \to C^\infty(W_-)$ as the composition of the covariant derivative (think gradient) and clifford multiplication.
An application of the index theorem gives
\[
\text{Index}(\partial\!\!\!/)=\frac{1}{24}\int_X p_1(X)\,.
\]

Thus,
\[
\text{Sign}(X)=8\text{Index}(\partial\!\!\!/)=8\left(\text{dim}_\C(\text{ker}(\partial\!\!\!/)-
\text{dim}_\C(\text{coker}(\partial\!\!\!/)\right)\,.
\]
Since the Dirac operator is quaterionic linear, if follows that the signature of any spin $4$-manifold is divisible by $16$. This result is known as Rochlin's theorem.
It also establishes that the trefoil knot is not slice, since the existence of a slice disk would imply the existence of a closed, spin $4$-manifold with signature negative eight.


We now return to the question of constructing the manifold $V_{\text{\rm res}}$. The starting point is to consider the singular algebraic variety
\[
V:= \{(x,y,z)\in \C^3 | x^2 + y^3 + z^5=0 \}\,.
\]
This has an isolated singularity. By definition the link of this singularity is the boundary of a neighborhood of the singularity, i.e.
\[
\Sigma = \Sigma(2,3,5) := \{(x,y,z)\in \C^3 | x^2 + y^3 + z^5=0 \ \text{and} \ \text{max}(|x|,|y|)=1 \}\,.
\]

The variety $V$ is just an open cone on $\Sigma(2,3,5)$, as one can see via the map
$[0,\infty) \times \Sigma(2,3,5) \to V$ by $(t,(x,y,z))\mapsto (t^{15}x, t^{10}y, t^{6}z)$. The closed cone on $\Sigma(2,3,5)$ is denoted by $\bar V$. It is just the set of points in $V$ with $\text{max}(|x|,|y|)\le 1$. The map
\[
\Sigma \to S^3:=\{(x,y)\in\C^2 | \text{max}(|x|,|y|)=1 \}\,,
\]
given by $(x,y,z)\mapsto (x,y)$ is a $5$-fold branched cover. The branch set is
$T(2,3):=\{(x,y)\in S^3 | x^2+y^3=0 \}$. The torus $\{(x,y)\in S^3 | |x|=|y|=1\}$ becomes a standard torus when $S^3-\{(0,1)\}$ is identified with $\R^3$. The branch set can be parameterized as $x=ie^{it/2}$, $y=e^{it/3}$, and this wraps twice around the torus in one direction while wrapping three times around the torus in the other direction. Thus, the branch set is a trefoil knot.


We will construct the manifold $V_{\text{\rm res}}$ by resolving the singularity in $V$. The manifold $\bar V_{\text{\rm res}}$ is just the resolution of the closed cone $\bar V$. Recall that a resolution of a singular algebraic variety $V$ is a smooth algebraic variety $V_{\text{\rm res}}$ together with an algebraic map to $V$ that is an isomorphism away from a singular set with codimension $1$ or more.
From the definition it is clear that $V$ can be expressed as a $5$-fold, or a $3$-fold, or a $2$-fold branched cover over $\C^2$. Following \cite{hkk} we will
use the $2$-fold branched cover description and take a sequence of blow-ups of
the base $\C^2$ to remove:
\begin{enumerate}
\item Singularities in the branch set,
\item Tangency between components of the branch set,
\item Triple point intersections of the branch set,
\item Intersections between different odd multiplicity components of the branch set.
\end{enumerate}
It will be easier to resolve the singularities in the pull-back of the cover to this multiple blow-up.

We begin this procedure by describing the blow-up process. Notice that the negative
degree line bundles over $\CP^1$ can be imbedded in the trivial rank two bundle:
\[
L_{-n}\to \{(x,y,[p:q])\in \C^2\times\CP^1 | \det \begin{bmatrix}x & y \\ p^n & q^n \end{bmatrix} = 0\}\,,
\]
for $n\ge 0$ given by sending the equivalence class $[x,y,\zeta]$ in $L_{-n}$ to $(\zeta x^n, \zeta y^n, [x:y])$. For $n=1$ the projection on the first two coordinates
is an algebraic isomorphism to $\C^2-\{(0,0)\}$ away from the exceptional divisor,
$E=\{(0,0,[p:q])\}\subseteq L_{-1}$.  This is called the blow-up at $(0,0)$. Any neighborhood of a smooth algebraic surface isomorphic to $\C^2$ can be replaced
by a copy of $L_{-1}$ and this is the blow-up.

\begin{figure}[!ht]  
\hskip20bp
\includegraphics[width=125mm]{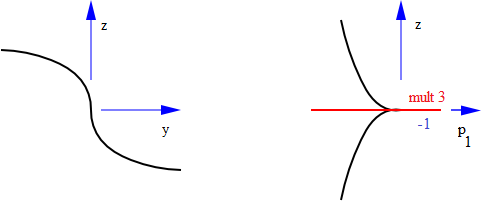}\caption{The branch set and inverse image}\label{b1b2}
\end{figure}

Our sequence of blow-ups begins with $U_0=\C^2$ and $B_0=\{(y,z)\in \C^2 | y^3+z^5=0 \}$, the branch set of the $2$-fold cover from $V$. This branch set has a singularity at $(0,0)$, so we blow up this point to get $U_1 = L_{-1}$. To see the new branch set, we will work in the $q_1:=q=1$ chart. Thus $y=p_1z$ and we rewrite everything using the $(p_1,z)$-coordinates. The equation $y^3+z^5=0$ becomes $(p_1^3+z^2)z^3=0$,
so we define
\[
B_1=\{(p_1,z) | (p_1^3+z^2)z^3=0 \}\,.
\]
To get the full copy of everything, we would also need to consider the other chart
(the $p_1=1$ chart) in the blow-up. As nothing significant occurs in this chart we do not write it. The set $B_1$ has two components corresponding to $p_1^3+z^2=0$ and $z^3=0$. The first component is non-compact, and the second component is a copy of the exceptional divisor in the blow-up with multiplicity $3$. In figure \ref{b1b2} below we display $B_0$ and $B_1$. We draw the compact components in red and the non-compact components in black. We further mark the multiplicity of each compact component in red next to the component.

We have also labeled each compact component with its self-intersection number in blue.
A generic section of the bundle $L_d$ intersects the zero section in exactly $d$ points counted with orientation. This is easy to see when $d\ge 0$ because a holomorphic section can be represented by a homogeneous polynomial of degree $d$, say $f(p,q)$ via $(s([p:q])=[p,q,f(p,q)]$. Such a polynomial has exactly $d$ zeros in the projective plane. For a generic polynomial we can work in the $q=1$ chart and this is just the count of the zeros of a degree $d$ polynomial. Using representatives of points in $L_{-n}$, say $[x,y,\zeta]$ with $|x|^2+|y|^2=1$ we can consider a section given by $(s([p:q])=[p,q,f(\bar p,\bar q)]$, where $f$ has degree $n$. Thus the zero section of $L_{-n}$ has self-intersection $-n$.

The next thing to do is to blow-up the point $(0,0,[0:1])$ in $U_1$. As the $q_1=1$ chart is isomorphic to $\C^2$ using $(p_1,z)$ coordinates, this amounts to replacing the $(p_1,z)$ patch with
\[
\{(p_1,z,[p_2:q_2])\in \C^2\times\CP^1 | \det \begin{bmatrix}p_1 & z \\ p_2 & q_2 \end{bmatrix} = 0\}\,.
\]
Working in the $p_2=1$ chart, we get $z=p_1q_2$, so the inverse image of the branch set $B_2$ is given by $p_1^5q_2^3(p_1+q_2^2)=0$. Once again we don't need to analyze the other chart. Two components of this set are tangent at $p_1=0$, $q_2=0$, thus we need to blow up this point to get $U_3$ by replacing the $(p_1,q_2)$ chart by
\[
\{(p_1,q_2,[p_3:q_3])\in \C^2\times\CP^1 | \det \begin{bmatrix}p_1 & q_2 \\ p_3 & q_3 \end{bmatrix} = 0\}\,.
\]

\begin{figure}[!ht]  
\hskip10bp
\includegraphics[width=125mm]{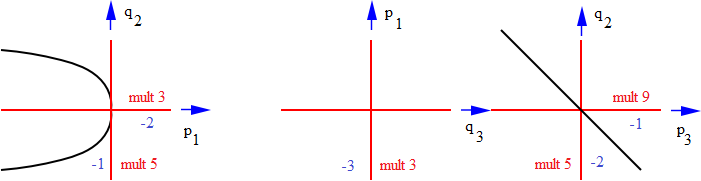}\caption{The branch set after two and three blow-ups}\label{b3b4}
\end{figure}

This time we need to analyze the $p_3=1$ chart and the $q_3=1$ chart. The $p_3=1$ chart contains the multiplicity $3$ cycle (compact component) as one sees by the formula $p_1^9q_3^3(1+p_1q_3^2)=0$ for $B_3$ in this chart. This does not appear in the $q_3=1$ chart, however the $q_3=1$ chart contains a triple point. The formula in this chart  is $p_3^5q_2^9(p_3+q_2)=0$. Figure \ref{b3b4} displays one chart intersecting $B_2$ on the left and two charts meeting $B_3$ on the right.

\begin{figure}[!ht]  
\hskip10bp
\includegraphics[width=125mm]{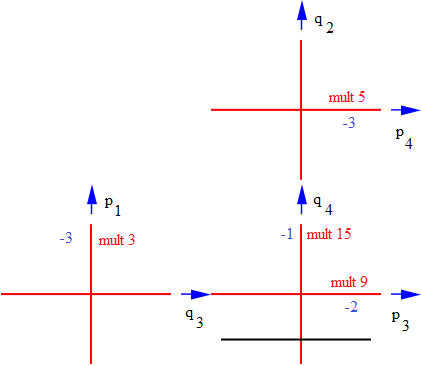}\caption{The branch set after four blow-ups}\label{b5}
\end{figure}

The next step is to blow up the triple point in the $q_3=1$ chart. This amounts to replacing the $(p_3,q_2)$ coordinate patch with
\[
\{(p_3,q_2,[p_4:q_4])\in \C^2\times\CP^1 | \det \begin{bmatrix}p_3 & q_2 \\ p_4 & q_4 \end{bmatrix} = 0\}\,.
\]
In the $p_4=1$ chart,  $B_4$ is given by $p_3^15q_4^9(1+q_4)=0$. We also keep track of the $q_4=1$ chart and the $(q_3,p_1)$ chart that was unchanged by the blow-up. The result is shown in figure \ref{b5}.

We continue the blow-up process everywhere two components with odd multiplicity intersect.
The local model for such an intersection is $\{(s,t)\in \C^2 | s^{2n+1}t^{2m+1}=0\}$. Blowing-up the intersection point will replace this chart with
\[
\{(s,t,[p:q])\in \C\times\CP^1 \det \begin{bmatrix}s & t \\ p & q \end{bmatrix} = 0\}\,.
\]
In the $p=1$ chart $t=sq$ and the branch set is described by \newline $s^{2n+2m+2}q^{2m+1}=0$. In the $q=1$ chart $s=tp$ and the branch set is described by $p^{2n+1}t^{2n+2m+2}=0$. Thus we see that blowing up a point of intersection between two cycles each with odd multiplicity results in a branch set in which the two cycles no longer intersect, but where a new cycle with multiplicity equal to the sum of the original multiplicities intersects the inverse image of each of the
original cycles transversely in one point.

\begin{figure}[!ht]  
\hskip10bp
\includegraphics[width=125mm]{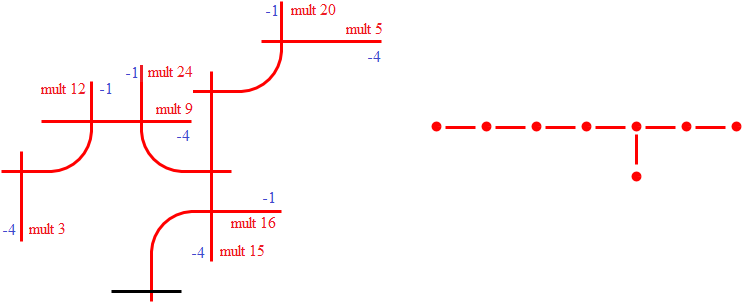}\caption{The branch set after all blow-ups}\label{b9}
\end{figure}

At this point we will no longer write out all of the coordinates and the coordinate charts. We will also combine all of the charts together into one picture.
Blowing up at the four points where odd multiplicity components intersect leads to the branch set in the spaces $U_5$, $U_6$, $U_7$, and $U_8$. The result in $U_8$ is displayed in figure \ref{b9}. This figure also includes a graph with one vertex for each cycle and two cycles connected by an edge if they intersect.

In general, the pull back of a space $E$ with a map $f:E\to Z$ by a map $g:Y\to Z$ is simply
\[
g^*E:=\{(y,e)\in Y\times E | g(y)=f(e)\}\,.
\]
The composition of all of the blow-up maps is a map $\pi:U_8\to U_0:=\C^2$. We now take the pull-back of the $2$-fold branched cover $f:V\to \C^2$ to obtain
\[
\pi^*V:=\{(u,(x,y,z))\in U_8\times V | \pi(u)=f(x,y,z)=(y,z) \}\,.
\]

This pull-back can be covered by charts in which it takes the form
\[
\{(x,s,t) | x^2 = s^{2n+1}t^{2m} \}\,.
\]

While these varieties are not smooth, it is not difficult to resolve the singularities here. Indeed, the map
\[
\C^2\to\{(x,s,t) | x^2 = s^{2n+1}t^{2m} \}\ \text{given by} \ (r,t)\mapsto (r^{2n+1}t^m,r^2,t)\,,
\]
does the job. Each point in the codomain with $t\neq 0$ has a unique inverse image, as does the point $(0,0,0)$. The remaining points ($t=0$ and $s\neq 0$) each have two inverse images. All of this takes place over the cycles with even multiplicity in the eight-fold blow-up $U_8$, and all of the blow-ups take place over the image of the original singular point in $V$ under the $2$-fold branched cover. (The first blow-up was over the image of the singular point and the remaining ones were on the exceptional divisors (cycles) that were introduced in the blow-up process.) Fitting copies of $\C^2$ together according to the transition maps identifying the charts in the pull-back $\pi^*V$ would be one way to describe the resolution $V_{\text{\rm res}}$. An alternate way to view this is to notice that the result is a $2$-fold branched cover over the even multiplicity cycles in $U_8$. Working in our chart the branched cover is just $(x,s,t)\mapsto(s,t)$. Translating to the resolution gives
$(r,t)\mapsto (r^2,t)$.

To analyze the intersection pairing on $V_{\text{\rm res}}$, we just need to analyze
the branched cover in neighborhoods of the cycles. There are two cases to consider -- cycles with odd multiplicity, and cycles with even multiplicity. Notice that the branch set (the union of even multiplicity components) intersects each odd multiplicity cycle in exactly two points and these cycles have self-intersection number $-4$. It follows that a local model of the branched cover in neighborhoods of these cycles is given by
\[
\begin{aligned}
\{(x,y,[p:q])\in &\C^2\times\CP^1 | \det \begin{bmatrix}x & y \\ p^2 & q^2 \end{bmatrix} = 0\} \\
&\to \{(z,w,[p:q])\in \C^2\times\CP^1 | \det \begin{bmatrix}x & y \\ p^4 & q^4 \end{bmatrix} = 0\}\,,
\end{aligned}
\]
with $(x,y,[p:q])\mapsto(x^2,y^2,[p:q])$.
It follows that the inverse image of each odd multiplicity cycle in $U_8$ is a cycle with self-intersection number $-2$ in
$V_{\text{\rm res}}$.


Now notice that each even multiplicity cycle in $U_8$ has self-intersection number $-1$ and is part of the branch locus of the $2$-fold branched cover  $V_{\text{\rm res}}\to U_8$. It follows that a local model of the cover over neighborhoods of these cycles is given by
\[
\begin{aligned}
\{(x,y,[p:q])\in &\C^2\times\CP^1  | \det \begin{bmatrix}x & y \\ p^2 & q^2 \end{bmatrix} = 0\} \\
&\to \{(x,y,[u:v])\in \C^2\times\CP^1 | \det \begin{bmatrix}x & y \\ u & v \end{bmatrix} = 0\}\,,
\end{aligned}
\]
with $(x,y,[p:q])\mapsto(x,y,[p^2:q^2])$.
It follows that the inverse image of each even multiplicity cycle in $U_8$ is a cycle with self-intersection number $-2$ in
$V_{\text{\rm res}}$.

Thus the graph of the cycles in $V_{\text{\rm res}}$ is exactly the same as the graph of the cycles in $U_8$ displayed in figure \ref{b9} except all of the self-intersection numbers are $-2$. Thus we see that the intersection form is the famous $E_8$ matrix:
\[
E_8=\begin{bmatrix}
-2 &1 & & & & & & \\
1& -2& 1 & & & & & \\
 & 1& -2 & 1 & & & & \\
  & & 1&-2&1& & & \\
  & & &1&-2&1& 0&1\\
  & & & &1&-2&1& \\
  & & & &0&1&-2&1\\
  & & & &1& & &-2
\end{bmatrix}
\]

\section*{Conclusion}
Michael Atiyah and Isadore Singer announced the Index theorem 50 years ago. This theorem connects deep topological properties of spaces to the structure of elliptic differential operators. Since the two different sides of the equation come from different branches of mathematics, the theorem serves as a bridge that connects different mathematical results. The index theorem has many applications and many generalizations, and is one of the great mathematical achievements of the last century. Atiyah and Singer were awarded the Abel Prize in 2004, in part, for their index theorem.

When I was an undergraduate, Ed Ihrig showed me the index theorem in an undergraduate geometry class. (Thank Ed.) It inspired me to learn more about topology. I hope this exposition inspires others to explore and learn more about the wonderful world of geometry.

\bibliography{ind}

\end{document}